\documentclass[10pt]{article}

\usepackage{a4wide}
\usepackage{amssymb}
\usepackage{amsfonts}
\usepackage{amsmath}
\input xy
\xyoption{arrow} \xyoption{matrix}

\date{}

\newtheorem{proposition}{Proposition}[section]
\newtheorem{theorem}[proposition]{Theorem}
\newtheorem{lemma}[proposition]{Lemma}

\newtheorem{corollary}[proposition]{Corollary}

\def\GK{{\rm  GK}\,}
\def\Kdim{{\rm K.dim }\,}

\def\der{\partial }

\def\nFM0{{\nu }_{F,M_0}}
\def\nFN0{{\nu }_{F,N_0}}
\def\nGN0{{\nu }_{G,N_0}}

\def\N0{ {\bf N}_0 }

\def\t{\otimes}
\def\g{\gamma}
\def\v{\varphi}
\def\ra{\rightarrow}

\def\lra{\leftrightarrow}
\def\Xpm{X^{\pm }}

\def\s{\sigma}
\def\Z{\mathbb{Z}}

\def\l1{{\lambda}_1}

\def\a{\alpha}
\def\a0{ {\alpha }_0}
\def\a1{ {\alpha }_1}

\def\l{\lambda}

\def\nFGM0{{\nu }_{F,G,M_0}}

\def\nFN0{{\nu}_{F,N_0}}

\def\sm{{\sigma}^m}

\def\sm1{{\sigma}^{-1}}

\def\smtp1{{\sigma}^{-t+1}}

\def\S1{S^{-1}}

\def\Xpm1{X^{\pm 1}_1}

\def\sPM1{{\sigma }^{\pm 1}}
\def\sMP1{{\sigma }^{\mp 1 }}

\def\d{\delta}

\def\di{{\rm d.ind}}

\def\L{\Lambda}

\def\CD{{\cal D}}

\def\Ytm1{Y^{t-1}}
\def\Yim1{Y^{i-1}}

\def\CM{{\cal M}}
\def\CN{{\cal N}}

\def\ass{{\rm ass}}

\def\supp{{\rm supp }}

\def\bK{\overline{K}}

\def\dim{{\rm dim }}
\def\char{{\rm char }}

\def\ker{ {\rm ker } }

\def\CJ{ {\cal J}}
\def\gr{ {\rm gr} }

\def\SL2Z{ {\rm SL}_2({\bf Z}) }

\def\Gp1{ G^{1 , 1 } }
\def\P11{ P^{-1 , 1 } }
\def\Pp1{ P^{1 , 1 } }

\def\nCLsr{{}^\nu\kern-2pt {\cal L}^{\sigma , \rho  }}
\def\nP{{}^\nu \kern-2pt P}
\def\nL{{}^\nu\kern-2pt L}
\def\nLL{{}^\nu\kern-2pt \Lambda}
\def\nPsr{{}^\nu\kern-2pt P^{\sigma , \rho  }}
\def\nLsr{{}^\nu\kern-2pt L^{\sigma , \rho  }}
\def\nuCL{{}^\nu\kern-2pt  {\cal L}}
\def\nCLsr{{}^\nu\kern-2pt {\cal L}^{\sigma , \rho  }}
\def\nCL1m{{}^\nu\kern-2pt {\cal L}^{-1 , 1  }}
\def\x1nu{x^\frac{1}{\nu}}
\def\xm1nu{x^{-\frac{1}{\nu}}}

\def\CN{{\cal N}}
\def\ra{\rightarrow }

\def\CB{{\cal B}}

\def\coker{{\rm coker}}

\def\CC{ {\cal C}}

\def\CP{ {\cal P}}

\def\nAM0{{\nu }_{{\cal A},M_0}}
\def\nAN0{{\nu }_{{\cal A},N_0}}

\def\Kdim{ {\rm Kdim } }

\def\End{ {\rm End }}

\def\CJ{ {\cal J }}

\def\CP{ {\cal P }}

\def\bI{\overline{I}}

\def\spec{{\rm spec}}

\def\bx{\overline{x}}
\def\by{\overline{y}}

\def\ga{\mathfrak{a}}
\def\gb{\mathfrak{b}}
\def\gc{\mathfrak{c}}
\def\gd{\mathfrak{d}}

\def\gn{\mathfrak{n}}
\def\gm{\mathfrak{m}}
\def\gp{\mathfrak{p}}
\def\gq{\mathfrak{q}}
\def\gr{\mathfrak{r}}

\def\bJ{\overline{J}}

\def\SL{{\rm SL}}

\def\Spec{{\rm Spec}}

\def\Ext{{\rm Ext}}

\def\di!{\frac{\der^i}{i!}}
\def\dik!{\frac{\der^k_i}{k!}}
\def\hA{\widehat{A}}

\def\id{{\rm id}}

\def\Max{{\rm Max}}

\def\N{\mathbb{N}}

\def\0{\overline{0}}
\def\1{\overline{1}}

\def\Ln1{\L_{n,\overline{1}}}

\def\a1{a_{\overline{1}}}

\def\S{\Sigma}

\def\vn1{\overrightarrow{n-1}}

\def\im{{\rm im}}

\def\mA{\mathbb{A}}

\def\soc{{\rm soc}}
\def\Sub{{\rm Sub}}
\def\SSub{{\rm SSub}}
\def\Inc{{\rm Inc}}
\def\Min{{\rm Min}}
\def\csupp{{\rm csupp}}

\def\F{\mathbb{F}}

\def\mS{\mathbb{S}}

\def\mI{\mathbb{I}}
\def\mrD{\mathrm{D}}
\def\ann{{\rm ann}}
\def\lann{{\rm l.ann}}
\def\rann{{\rm r.ann}}
\def\Cen{{\rm Cen}}
\def\clKdim{{\rm cl.Kdim}}
\def\pd{{\rm pd}}
\def\hmS{\widehat{\mathbb{S}}}
\def\hht{{\rm ht}}
\def\cht{{\rm cht}}
\def\lgldim{{\rm l.gldim}}
\def\rgldim{{\rm r.gldim}}

\def\bM{\overline{M}}

\def\wdim{{\rm wdim}}

\begin{document}

\author{V. V. \  Bavula 
}

\title{The algebra of one-sided inverses of a polynomial algebra}


\maketitle
\begin{abstract}
We study in detail the 
algebra $\mS_n$ in the title which is  an algebra obtained from a
polynomial algebra $P_n$ in $n$ variables by adding commuting,
{\em left} (but not two-sided) inverses of the canonical
generators of $P_n$. The algebra $\mS_n$ is non-commutative and
neither left nor right Noetherian but the set of its ideals
satisfies the a.c.c., and the ideals {\em commute}. It is proved
that the classical Krull dimension of $\mS_n$ is $2n$; but the
weak and the global dimensions of $\mS_n$ are $n$. The prime and
maximal spectra of $\mS_n$ are found, and the simple
$\mS_n$-modules are classified. It is proved that the algebra
$\mS_n$ is central, prime, and {\em catenary}. The set $\mI_n$ of
idempotent ideals of $\mS_n$ is found explicitly. The set $\mI_n$
is a finite distributive lattice and the number of elements in the
set $\mI_n$ is equal to the {\em Dedekind} number $\gd_n$.



$${\bf Contents}$$
\begin{enumerate}
\item Introduction.

\item The ideals of $\mS_n$ commute and satisfy the a.c.c..

\item Classification of simple $\mS_n$-modules.

\item The prime and maximal spectra of the algebra $\mS_n$.

 \item Left or right Noetherian factor algebras of $\mS_n$.

\item The weak and  global dimensions of the algebra $\mS_n$.

\item Idempotent ideals of the algebra $\mS_n$.

\end{enumerate}
\end{abstract}


\section{Introduction}
Throughout, ring means an associative ring with $1$; module means
a left module;
 $\N :=\{0, 1, \ldots \}$ is the set of natural numbers; $K$ is a
field 
and  $K^*$ is its group of units;
$P_n:= K[x_1, \ldots , x_n]$ is a polynomial algebra over $K$;
$\der_1:=\frac{\der}{\der x_1}, \ldots , \der_n:=\frac{\der}{\der
x_n}$ are the partial derivatives ($K$-linear derivations) of
$P_n$.

$\noindent $

The 
 algebras $\mS_n$ (see the definition below) appear
naturally when one wants to develop a {\em theory of one-sided
localizations}. Let me give an example. Let $K[x]$ be a polynomial
algebra in a variable $x$ over the field $K$. When we invert the
element $x$ the resulting algebra $K[x,x^{-1}]$ has the same
properties as $K[x]$. This is not the case when we invert the
element $x$ on {\em one side} only, say, on the left: $yx=1$. Then
the algebra $\mS_1:= K\langle x,y \, | \, yx=1\rangle$ has very
different properties from the polynomial algebra $K[x]$. It is
non-commutative, not left and right Noetherian, not a domain, it
contains the ring of infinite dimensional matrices, etc. Moreover,
the algebra $\mS_1$ has properties that are a mixture of the
properties of the three algebras: $K[x_1]$, $K[x_1,x_2]$, and the
{\em Weyl algebra} $A_1:=K\langle x, \der \, | \, \der x -x\der
=1\rangle$ if $\char (K)=0$ (for example, as it is proved in the
paper, the ideals of the algebra $\mS_1$ commute, and each proper
ideal of $\mS_1$ but one is a unique product of maximal ideals
(counted with multiplicity), and the lattice of ideals is
distributive; the classical Krull dimension of $\mS_1$ is 2; the
global homological dimension of $\mS_1$ is 1; the Gelfand-Kirillov
dimension of $\mS_1$ is 2). The algebra $\mS_1$ is a well-known
primitive algebra \cite{Jacobson-StrRing}, p. 35, Example 2. Over
the field
 $\mathbb{C}$ of complex numbers, the completion of the algebra
 $\mS_1$ is the {\em Toeplitz algebra} which is the
 $\mathbb{C}^*$-algebra generated by a unilateral shift on the
 Hilbert space $l^2(\N )$ (note that $y_1=x_1^*$). The Toeplitz
 algebra is the universal $\mathbb{C}^*$-algebra generated by a
 proper isometry.

$\noindent $

{\it Definition}. The 
algebra $\mathbb{S}_n$ of one-sided inverses of $P_n$ is an
algebra generated over a field $K$ 
by $2n$
elements $x_1, \ldots , x_n, y_n, \ldots , y_n$ that satisfy the
defining relations:
$$ y_1x_1=\cdots = y_nx_n=1 , \;\; [x_i, y_j]=[x_i, x_j]= [y_i,y_j]=0
\;\; {\rm for\; all}\; i\neq j,$$ where $[a,b]:= ab-ba$, the
commutator of elements $a$ and $b$.

$\noindent $

By the very definition, the algebra $\mS_n$ is obtained from the
polynomial algebra $P_n$ by adding commuting, left (or right)
inverses of its canonical generators.
 Clearly,
$\mathbb{S}_n=\mS_1(1)\t \cdots \t\mS_1(n)\simeq \mathbb{S}_1^{\t
n}$ where $\mS_1(i):=K\langle x_i,y_i \, | \, y_ix_i=1\rangle
\simeq \mS_1$ and $$\mS_n=\bigoplus_{\alpha , \beta \in \N^n}
Kx^\alpha y^\beta$$ where $x^\alpha := x_1^{\alpha_1} \cdots
x_n^{\alpha_n}$, $\alpha = (\alpha_1, \ldots , \alpha_n)$,
$y^\beta := y_1^{\beta_1} \cdots y_n^{\beta_n}$, $\beta =
(\beta_1,\ldots , \beta_n)$. In particular, the algebra $\mS_n$
contains two polynomial subalgebras $P_n$ and $Q_n:= K[y_1, \ldots
, y_n]$ and is equal,  as a vector space,  to their tensor product
$P_n\t Q_n$.  The canonical generators $x_i$, $y_j$ $(1\leq
i,j\leq n)$ determine the ascending filtration $\{ \mS_{n, \leq
i}\}_{i\in \N}$ on the algebra $\mS_n$ in the obvious way (i.e. by
the total degree of the generators): $\mS_{n, \leq i}:=
\bigoplus_{|\alpha |+|\beta |\leq i} Kx^\alpha y^\beta$ where
$|\alpha | \; = \alpha_1+\cdots + \alpha_n$ ($\mS_{n, \leq
i}\mS_{n, \leq j}\subseteq \mS_{n, \leq i+j}$ for all $i,j\geq
0$). Then $\dim (\mS_{n,\leq i})={i+2n \choose 2n}$ for $i\geq 0$,
and so the Gelfand-Kirillov dimension $\GK (\mS_n )$ of the
algebra $\mS_n$ is equal to $2n$. It is not difficult to show
(Lemma \ref{c26Nov8}) that the algebra $\mS_n$ is neither left nor
right Noetherian. Moreover, it contains infinite direct sums of
left and right ideals.

$\noindent $

 Another (left and right) non-Noetherian  algebras, so-called, the
 Jacobian algebras $\mA_n$ (introduced in \cite{Bav-Jacalg}), appear as a localization {\em not} in the sense of
 Ore of the Weyl algebras $A_n$. The general construction, proposed by the author,  is as
 follows: given an algebra $A$, an $A$-module $M$, and a set $S=\{
 a_i\}_{i\in I}$ of elements of the algebra $A$ such that the maps
 $a_{i, M}: M\ra M$,  $ m\mapsto a_im$, are invertible. The
 subalgebra $S^{-1}_M*A$  of $\End_K(M)$ generated by the image of the
 algebra $A$ in $\End_K(M)$ and the elements $\{a_{i,
 M}^{-1}\}_{i\in I}$ can be seen as a new way of localizing  the algebra $A$. Clearly,
 $S^{-1}_M*(S^{-1}_M*A)=S^{-1}_M*A$, and the factor algebra $A/\ann_A(M)$ of $A$
 modulo the annihilator $\ann_A(M)$ of the $A$-module $M$ is a
 subalgebra of $S^{-1}_M*A$. In general, as the example of the
 Jacobian algebras $\mA_n$ shows \cite{Bav-Jacalg},  the algebras $A$ and $S^{-1}_M*A$ have
 different properties, and the localized algebra $S^{-1}_M*A$ is not a left
 or right localization of the algebra $A$ in the sense of Ore.

$\noindent $

{\it Definition}, \cite{Bav-Jacalg}. When  $\char (K)=0$,  the
{\em Jacobian algebra} $\mA_n$ is the subalgebra of $\End_K(P_n)$
generated by the {\em Weyl algebra} $A_n:=K\langle x_1, \ldots ,
x_n, \der_1, \ldots , \der_n\rangle$ and the elements $H_1^{-1},
\ldots , H_n^{-1}\in \End_K(P_n)$ where $H_1:= \der_1x_1, \ldots ,
H_n:= \der_nx_n$.

$\noindent $

Clearly, $\mA_n = \mA_1(1) \t \cdots \t \mA_1(n) \simeq \mA_1^{\t
n }$ where $\mA_1(i) := K\langle x_i, \der_i , H_i^{-1}\rangle
\simeq \mA_1$. The algebra $\mA_n$ contains all the  integrations
$\int_i: P_n\ra P_n$, $ p\mapsto \int p \, dx_i$, since  $\int_i=
x_iH_i^{-1}$. In particular, the algebra $\mA_n$ contains all
(formal) integro-differential operators with polynomial
coefficients. The Jacobian algebra $\mA_n$ appeared in my study of
the group of polynomial automorphisms and the Jacobian Conjecture,
which is a conjecture that makes sense {\em only} for polynomial
algebras in the class of all commutative algebras
\cite{Bav-inform}. In order to solve the Jacobian Conjecture,   it
is reasonable to believe that one should create a technique which
makes sense {\em only} for polynomials;  the Jacobian algebras are
a step in this direction (they exist for polynomials but make no
sense even for Laurent polynomials). The Jacobian algebras were
studied in detail in \cite{Bav-Jacalg}. Their relevance to 
 $\mS_n$ is the fact that $\mS_n$ is a subalgebra of $\mA_n$ (Lemma
\ref{a26Nov8}), and this fact makes it possible to shorten proofs
of several results on $\mS_n$. Moreover, there are many parallels
between these two classes of algebras.

$\noindent $

Let us describe main results of this paper.

\begin{itemize}
 \item  (Proposition \ref{b21Dec8}, Theorem \ref{23Dec8}) {\em The
 algebra $\mS_n$ is central, prime, and catenary. Every nonzero
 ideal  of $\mS_n$ is an essential left and right submodule of}
 $\mS_n$.
\item (Proposition \ref{c21Dec8}). {\em Let $A$ be a $K$-algebra.
 Then the algebra $\mS_n\t A$ is a prime algebra iff the algebra
$A$ is a prime algebra.}
 \item (Theorems \ref{8Dec8} and \ref{b23Dec8}). {\em The ideals of
 $\mS_n$ commute ($IJ=JI$),  and the set of ideals of $\mS_n$ satisfy the a.c.c..}
 \item (Theorem \ref{f21Dec8}) {\em The classical Krull dimension $\clKdim (\mS_n)$ of $\mS_n$ is $2n$.}
 \item (Theorem \ref{25Feb9}) {\em For each ideal $\ga$ of
  $\mS_n$, the set $\Min (\ga )$ of the minimal primes over $\ga$
 is a finite, non-empty set.}
 \item (Theorem \ref{22Dec8}) $\hht (\gp ) +\cht (\gp) = \clKdim
 (\mS_n )$, {\em for all prime ideals $\gp$ of $\mS_n$. Formulae for the
 height $\hht (\gp )$ and the co-height $\cht (\gp )$ are found
 explicitly (via combinatorial data).}
 \item (Theorems \ref{A11Dec8}, Corollary  \ref{26Dec8})  {\em The
  weak homological dimension and the left and right
 global dimensions of $\mS_n$ are equal to $n$}.
 \item (Theorems \ref{21Dec8} and \ref{d21Dec8}) {\em The prime
 and the maximal spectra of $\mS_n$ are found.}
 \item (Theorem \ref{17Dec8}) {\em The simple $\mS_n$-modules are
 classified.}
 \item (Corollary \ref{c17Dec8}) {\em $\GK (M) \leq n$, for all simple
 $\mS_n$-modules $M$. Moreover, $\GK (M) \in \{0,1, \ldots , n\}$.}
 \item (Corollary \ref{b17Dec8}) {\em The annihilators separate
 the simple $\mS_n$-modules.}
 \item (Corollary \ref{a19Dec8}) {\em The polynomial algebra $P_n$
 is the only faithful, simple $\mS_n$-module.}
 \item (Corollary \ref{x28Dec8}) {\em $\GK (M)+\pd (M)=\lgldim (\mS_n)$ for all simple
 $\mS_n$-modules $M$.}
\item (Theorem \ref{28Dec8}) {\em The projective dimension $\pd
(M)$ of each the simple $\mS_n$-module $M$  is found explicitly.}
\item (Theorem \ref{17Dec8}.(4), Corollary \ref{d17Dec8}) {\em For
each simple $\mS_n$-module $M$, the endomorphism division algebra
$\End_{\mS_n}(M)$ is a finite field extension of $K$ and its
dimension over $K$ is equal to the multiplicity $e(M)$ of the
$\mS_n$-module $M$.} \item (Lemma \ref{a21Dec8}) {\em There are
precisely $n$ height one prime ideals of the algebra $\mS_n$ (they
are given explicitly).} \item (Theorem \ref{24Dec8}) {\em Let $I$
be an ideal of $\mS_n$. Then the factor algebra $\mS_n / I$ is
left (or right) Noetherian iff the ideal $I$ contains all the
height one primes of $\mS_n$.} \item (Corollary \ref{b9Jan9}) {\em
A prime ideal $\gp\neq \mS_n$  of $\mS_n$ is an idempotent ideal
($\gp^2=\gp$) iff  $\gp$ is contained  in all the maximal ideals
of $\mS_n$ iff $\gp$ is a sum of height one prime ideals of}
$\mS_n$. \item (Theorem \ref{9Jan9}, Corollary \ref{28Ma7}) {\em
The set $\mI_n$ of idempotent ideals of $\mS_n$ is found. The set
$\mI_n$  is a finite distributive lattice and the number of
elements in the set $\mI_n$ is equal to the Dedekind number
$\gd_n$.} \item (Theorem \ref{27Ma7}) {\em Each idempotent ideal
$\ga$ of $\mS_n$ distinct from $\mS_n$ is a unique product and a
unique intersection of the minimal (necessarily idempotent) prime
ideals over $\ga$.}
\end{itemize}
These results show that the algebra $\mS_n$ has properties that
are a mixture of properties of the polynomial algebra $P_{2n}$ and
the Weyl algebra $A_n$. This is not so  surprising when we look at
the defining relations of the algebras  $\mS_n$, $A_n$, and
$P_{2n}$.

The algebras $\mS_n$ are fundamental non-Noetherian algebras; they
are universal non-Noetherian algebras of their own kind in a
similar way as the polynomial algebras are universal in the class
of all the commutative algebras and the Weyl algebras are
universal in the class of algebras  of differential operators.

The algebra $\mS_n$ often appears as a subalgebra or a factor
algebra  of many non-Noetherian algebras. For example, $\mS_1$ is
a factor algebra of certain non-Noetherian down-up algebras as was
shown by Jordan \cite{Jordan-Down-up}  (see also \cite{Benk-Roby},
\cite{Kirk-Mus-Kuz}, \cite {Kirk-Kuz}).


\section{The ideals of $\mS_n$ commute and satisfy the a.c.c.}

In this section, it is proved that the algebra $\mS_n$ is neither
left nor right Noetherian (it contains infinite direct sums of
left and right nonzero ideals) but the set of all the ideals of
the algebra $\mS_n$ satisfies the a.c.c. (Theorem \ref{8Dec8}),
the ideals of the algebra $\mS_n$ commute (Theorem \ref{b23Dec8}).

The polynomial algebra $P_n$ is a left $\End_K(P_n)$-module, we
denote the action of a linear map $a\in \End_K(P_n)$ on an element
$p\in P_n$ either by $a (p)$ or by $a*p$. By the very definition,
the Jacobian algebra $\mA_n$ is a subalgebra of $\End_K(P_n)$.

\begin{lemma}\label{a26Nov8}
The algebra homomorphism $\mS_n \ra \mA_n$, $x_i\mapsto x_i$ , $
y_i\mapsto H_i^{-1}\der_i$, is a monomorphism when $\char (K)=0$.
\end{lemma}

{\it Proof}. In view of the natural isomorphisms $\mS_n\simeq
\mS_1^{\t n }$ and $\mA_n\simeq \mA_1^{\t n }$, it suffices to
prove the lemma when $n=1$ (we drop the subscript `1' in this case
here and everywhere if this does not lead to confusion). The
homomorphism is correctly defined since $H^{-1} \der x =1$. It
remains to show that its kernel is zero. Note that for each
natural numbers $i$ and $j$, we have
$$ x^iy^j\mapsto x^i(H^{-1} \der )^j = x^i \frac{1}{H(H+1)\cdots
(H+j-1)}\der^j.$$ If an element $a= a_jy^j + a_{j+1}y^{j+1}+\cdots
+ a_{j+k} y^{j+k}\in \mS_1$ (where all $a_s\in K[x]$, and $k\geq
0$) belongs to the kernel then
$$ 0= a* x^j = a_jy^j * x^j =  a_j,$$
and so $a=0$, i.e. the kernel is zero. $\Box $

$\noindent $

We identify the algebra $\mS_n$ with its isomorphic copy in the
algebra $\mA_n$ via the above monomorphism. Then, $\mS_n\subset
\mA_n\subset \End_K(P_n)$, and so $P_n$ is a {\em faithful}
$\mS_n$-module when $\char (K)=0$.

\begin{corollary}\label{b26Nov8}
The $\mS_n$-module $P_n$ is simple and faithful.
\end{corollary}

{\it Proof}. We have to prove that the $\mS_n$-module $P_n$ is
simple. For $n=1$ and natural numbers $i$ and $j$:
\begin{equation}\label{yixj}
y^i*x^j= \begin{cases}
0& \text{if $j<i$},\\
x^{j-i}& \text{if $j\geq i$}.\\
\end{cases}
\end{equation}
If $p= \sum \l_\alpha x^\alpha\in P_n$ (where $\l_\alpha \in K$)
is a nonzero polynomial of degree, say $d$, then $\l_\beta \neq 0$
for some element $\beta \in \N^n$ such that $|\beta | = d$. Then
$\l_\beta^{-1} y^\beta *p=1$. This means that the $\mS_n$-module
$P_n$ is simple.

Suppose that $a*P_n=0$ for a nonzero element $a=\sum a_\alpha
y^\alpha\in A$ where $a_\alpha \in P_n$, we seek a contradiction.
Fix $\alpha$ such that $a_\alpha \neq 0$ and $|\alpha |$ is the
least possible. Then $0\neq a_\alpha = a*x^\alpha =0$, a
contradiction. Therefore, $P_n$ is a faithful $\mS_n$-module.
$\Box $

$\noindent $

Later, we will see that $P_n$ is the {\em only} simple and
faithful $\mS_n$-module (Corollary \ref{a19Dec8}).

$\noindent $

{\it Example}. Consider a vector space $V= \bigoplus_{i\in
\N}Ke_i$ and two shift operators on $V$, $X: e_i\mapsto e_{i+1}$
and $Y:e_i\mapsto e_{i-1}$ for all $i\geq 0$ where $e_{-1}:=0$. By
Corollary \ref{b26Nov8} and (\ref{yixj}), the subalgebra of
$\End_K(V)$ generated by the operators $X$ and $Y$ is isomorphic
to the algebra $\mS_1$ $(X\mapsto x$, $Y\mapsto y)$. By taking the
$n$'th tensor power $V^{\t n }=\bigoplus_{\alpha \in
\N^n}Ke_\alpha$ of $V$ we see that the algebra $\mS_n$ is
isomorphic to the subalgebra of $\End_K(V^{\t n })$ generated by
the $2n$ shifts $X_1, Y_1, \ldots , X_n, Y_n$ that act in
different directions.

$\noindent $

When $n=1$, by (\ref{yixj}), for each natural number $i$, the
product $x^iy^i$ is the projection onto the ideal $(x^i)$ of the
polynomial algebra $K[x]$ in the decomposition $K[x] =
(\bigoplus_{j<i} Kx^j)\bigoplus (x^i)$. Therefore, the elements of
the
algebra $\End_K(P_n)$ (where $i,j\in \N$): 
\begin{equation}\label{Eija}
E_{ij}:= \begin{cases}
x^{i-j}(x^jy^j-x^{j+1}y^{j+1})& \text{if $i\geq j$ },\\
y^{j-i}(x^jy^j-x^{j+1}y^{j+1})& \text{if $i<j$ },\\
\end{cases}
\end{equation}
are the matrix units, i.e. $E_{ij}*x^k = \d_{jk} x^i$, $k\geq 0$,
where $\d_{jk}$ is the Kronecker delta. In particular, $E_{ij}
E_{kl}= \d_{jk}E_{il}$ for all natural numbers $i$, $j$, $k$,  and
$l$. Therefore, the subring (without 1) $$F:= \bigoplus_{i,j\in
\N}KE_{ij}$$ of $\mS_1$ is canonically isomorphic to the ring
(without 1) of infinite dimensional matrices $$M_\infty (K) :=
\varinjlim M_d(K)=\bigoplus_{i,j\in \N}KE_{ij}$$ (via $F\ra
M_\infty (K)$, $E_{ij}\mapsto E_{ij}$) where $E_{ij}$ are the
matrix units of $M_\infty (K)$ and $M_d(K):= \bigoplus_{i,j=1}^n
KE_{ij}$ is the ring of $d$-dimensional matrices over $K$.

We have another presentation of the matrix units: 
\begin{equation}\label{Eijb}
E_{ij}:= \begin{cases}
(x^iy^i-x^{i+1}y^{i+1})x^{i-j}& \text{if $i\geq j$ },\\
(x^iy^i-x^{i+1}y^{i+1})y^{j-i}& \text{if $i<j$ }. \\
\end{cases}
\end{equation}
 The formula (\ref{Eijb}) can be verified directly
using the inclusion $\mS_1\subset \End_K(P_1)$. Now, combining
(\ref{Eija}) and (\ref{Eijb}) we can write 
\begin{equation}\label{Eijc}
E_{ij}= x^iy^j-x^{i+1}y^{j+1}=x^iE_{00}y^j, \;\; i,j\geq 0.
\end{equation}


{\bf The involution $\eta$ on $\mS_n$}. The algebra $\mS_n$ admits
the {\em involution}
$$ \eta : \mS_n\ra \mS_n, \;\; x_i\mapsto y_i, \;\; y_i\mapsto
x_i, \;\; i=1, \ldots , n,$$ i.e. it is a $K$-algebra
anti-isomorphism ($\eta (ab) = \eta (b) \eta (a)$ for all $a,b\in
\mS_n$) such that $\eta^2 = \id_{\mS_n}$, the identity map on
$\mS_n$. So, the algebra $\mS_n$ is {\em self-dual} (i.e. it is
isomorphic to its opposite algebra, $\eta : \mS_n\simeq
\mS_n^{op}$). This means that left and right algebraic properties
of the algebra $\mS_n$ are the same.

For $n=1$ and all $i,j\in \N$, 
\begin{equation}\label{eEij}
\eta (E_{ij} )=E_{ji}.
\end{equation}
This follows from (\ref{Eijc}). The involution $\eta$ acts on the
ring $F=M_\infty (K)$ as the transposition.  In general case,
\begin{equation}\label{eFF}
\eta (F_n) = F_n,
\end{equation}
where $F_n:= F^{\t n }=\bigotimes_{i=1}^n F(i)=\bigoplus_{\alpha ,
\beta \in \N^n}KE_{\alpha \beta}$, $F(i):= \bigoplus_{s,t\in
\N}KE_{st}(i)$,  $E_{\alpha \beta}:=\bigotimes_{i=1}^n E_{\alpha_i
\beta_i}(i)$ where $E_{\alpha_i
\beta_i}(i):=x_i^{\alpha_i}(1-x_iy_i)y_i^{\beta_i}\in F(i)\subset
\mS_1(i)$. Clearly, $E_{\alpha \beta}E_{\g \rho}=\d_{\beta, \g
}E_{\alpha \rho}$ for all elements $\alpha, \beta , \g , \rho \in
\N^n$ where $\d_{\beta , \g }$ is the Kronecker delta function.
The involution $\eta$ acts on the `matrix' ring $F_n$ as the
transposition: 
\begin{equation}\label{eEij1}
\eta (E_{\alpha \beta} )=E_{\beta \alpha}.
\end{equation}

{\bf The algebra $\mS_n$ is neither left nor right Noetherian}. By
(\ref{Eija}) and (\ref{Eijb}), for all $i,j\geq 0$,
\begin{equation}\label{xyEij}
xE_{ij}=E_{i+1, j}, \;\; yE_{ij} = E_{i-1, j}\;\;\; (E_{-1,j}:=0),
\end{equation}
\begin{equation}\label{xyEij1}
E_{ij}x=E_{i, j-1}, \;\; E_{ij}y = E_{i, j+1} \;\;\;
(E_{i,-1}:=0).
\end{equation}
 By
(\ref{xyEij}) and (\ref{xyEij1}), $F$ is an ideal of the algebra
$\mS_1$. Note that $F$ is an ideal of the algebra $\mA_1$
(Corollary 2.6,
 \cite{Bav-Jacalg})  (from this fact it also follows that  $F$ is an ideal of  the algebra
 $\mS_1$ since $F\subset \mS_1\subset \mA_1$).

 By
(\ref{xyEij}) and (\ref{xyEij1}), for all $i,j\geq 0$,
\begin{equation}\label{xyEij2}
xE_{ij}=E_{i+1, j+1}x, \;\; E_{ij}y=yE_{i+1, j+1}.
\end{equation}
By (\ref{xyEij}), for each $i\geq 0$ and $j\geq 0$, the left
$\mS_1$-module $\mS_1E_{ij}= \bigoplus_{k\geq 0} KE_{kj}$ is
isomorphic to the left $\mS_1$-module $K[x]$ via the isomorphism
$\mS_1E_{ij} \ra K[x]$, $E_{0j}\mapsto 1$ (and so $E_{kj}\mapsto
x^k$, $k\geq 0$).

\begin{lemma}\label{c26Nov8}
The algebra $\mS_n$ is neither left nor right Noetherian.
Moreover, it contains  infinite direct sums of left and right
ideals.
\end{lemma}

{\it Proof}. The algebra is self-dual, so it suffices to prove,
say,  the first statement of the lemma. Since $\mS_n\simeq \mS_1\t
\mS_{n-1}$, it suffices to prove the lemma when $n=1$. In this
case, $F= \oplus_{i\in \N} \mS_1E_{ii}\simeq \oplus_{i\in \N} P_1$
is the direct sum of infinitely many copies of the simple
$\mS_1$-module $P_1$. $\Box $

$\noindent $

{\bf The elements $x_i$ and $y_i$ of $\mS_n$
 are not regular}. Let $r$ be an element of a ring $R$. The
 element $r$ is called {\em regular} if $ \lann_R(r)=0$ and
 $\rann_r(r)=0$ where $ \lann_R(r):= \{ s\in R\, | \, sr=0\}$ is
 the {\em left annihilator} of $r$ and $ \rann_R(r):= \{ s\in R\, | \, rs=0\}$ is
 the {\em right annihilator} of $r$.

The next lemma shows that the elements $x$ and $y$ of the algebra
$\mS_1$ are not regular.
\begin{lemma}\label{a7Dec8}
\begin{enumerate}
\item $\lann_{\mS_1}(x) = \mS_1E_{00} = \bigoplus_{i\geq 0}
KE_{i,0}=\bigoplus_{i\geq 0} Kx^i(1-xy)$ and $\rann_{\mS_1}(x)=0$.
\item $\rann_{\mS_1}(y) = E_{00}\mS_1 =  \bigoplus_{i\geq 0}
KE_{0, i}=\bigoplus_{i\geq 0} K (1-xy)y^i$ and
$\lann_{\mS_1}(y)=0$.
\end{enumerate}
\end{lemma}

{\it Proof}. 1. $yx=1$ implies $\rann_{\mS_1}(x)=0$. Since $E_{00}
x= (1-xy)x= x-x=0$,
$$\lann_{\mS_1}(x) \supseteq \mS_1E_{00} = \bigoplus_{i\geq 0}
KE_{i,0}= \bigoplus_{i\geq 0} Kx^i(1-xy).$$ To prove the reverse
inclusion note that the right $K[x]$-module $\mS_1/ K[x]$ is a
direct sum of its right  submodules $M_i:= \bigoplus_{j=1}^\infty
Kx^iy^j+K[x]$, $ i\in \N$. It follows from the equalities $x^iy^j
x= x^iy^{j-1}$ that the kernel of the linear map in $\mS_1/ K[x]$
given by the multiplication by the element $x$ on the right is
equal to $\bigoplus_{i\geq 0} Kx^iy+K[x]$. Clearly, $\lann_{\mS_1}
(x) \subseteq \bigoplus_{i\geq 0} Kx^iy+K[x]=Ky\bigoplus
\bigoplus_{i\geq 0} KE_{i,0}\bigoplus K[x]$. Now, one can easily
find that $\lann_{\mS_1}(x) = \bigoplus_{i\geq 0} KE_{i,0}$, as
required.

2. The second statement follows from the first by using the
involution $\eta$:
\begin{eqnarray*}
 \rann_{\mS_1}(y)& =& \eta^2(\rann_{\mS_1}(y))= \eta (\lann_{\mS_1}(\eta (y)))= \eta (\lann_{\mS_1}(x)) =
 \eta (\mS_1 E_{00}) \\
 &= & \eta (E_{00})\mS_1 = E_{00}\mS_1=\bigoplus_{i\geq 0} K
 (1-xy)y^i,
\end{eqnarray*}
where we have used the fact that $\eta (E_{00})=E_{00}$ (see
(\ref{eEij})). $\Box $

It follows from (\ref{Eijc}) that 
\begin{equation}\label{mS1d}
\mS_1= K\bigoplus xK[x]\bigoplus yK[y]\bigoplus F,
\end{equation}
the direct sum of vector spaces. Then 
\begin{equation}\label{mS1d1}
\mS_1/F\simeq K[x,x^{-1}], \;\; x\mapsto x, \;\; y \mapsto x^{-1},
\end{equation}
since $yx=1$, $xy=1-E_{00}$ and $E_{00}\in F$.

$\noindent $

{\it Example}. Let $V=\bigoplus_{i\in \N}Ke_i$ be a vector space.
By taking the matrix of a linear map in $V$ with respect to the
basis $ \{ e_i\}_{i\in \N}$ for the vector space $V$, we identify
the algebra $\End_K(V)$ with the algebra of infinite matrices $\{
\sum_{i,j\in \N}a_{ij}E_{ij}$ is an infinite sum where for each
$j$ almost all scalars $a_{ij}=0\}$ where $E_{ij}$ are the matrix
units. Let $\gn := \sum_{i\geq 0}E_{i+1,i}$ and $\gm :=
\sum_{i\geq 1} E_{i-1,i}$. A matrix of the form $\sum_{i>
0}\l_{-i}\gm^i+\l_0E+\sum_{i> 0} \l_i\gn^i$ is called a {\em
multidiagonal matrix} where $\l_i\in K$ (and the sums are finite).
A matrix is called an {\em almost multidiagonal} if it becomes a
multidiagonal matrix by adding a finite sum $\sum \mu_{ij}E_{ij}$,
$\mu_{ij}\in K$. The set of all almost multidiagonal matrices is a
subalgebra of $\End_K(V)$ which is isomorphic to the algebra
$\mS_1$, by (\ref{mS1d}) and (\ref{yixj}) $(\gn \lra x$, $\gm \lra
y$, $E_{ij}\lra E_{ij})$.

$\noindent $

For an element $a$ of an algebra $A$, the subalgebra of $A$,
$\Cen_A(a):= \{ b\in A\, | \, ba=ab\}$, is called the {\em
centralizer} of the element $a$ in the algebra $A$. By
(\ref{xyEij}), (\ref{xyEij1}), and (\ref{mS1d}), 
\begin{equation}\label{Cenxy}
\Cen_{\mS_1}(x) = K[x]\;\; {\rm and}\;\; \Cen_{\mS_1}(y)=K[y].
\end{equation}

We say that a $K$-algebra $A$ is {\em central} if its centre
$Z(A)$ is $K$. We denote by $\CJ (A)$ the set of all the ideals of
the algebra $A$. An ideal $I$ of the algebra $A$ is called a {\em
proper} ideal if $I\neq 0, A$.  The {\em classical Krull
dimension} of the algebra $A$ is denoted  $\clKdim (A)$. $\spec
(A)$ and $\Max (A)$ denote the {\em prime spectrum} and the {\em
maximal spectrum} of the algebra $A$ respectively. A nonzero
polynomial $a\in K[x]$ is a {\em monic} polynomial if its leading
coefficient is 1. The {\em socle} $\soc (M)$ of a module $M$ is
the sum of its semi-simple submodules, if they exist, and is zero
otherwise.

\begin{proposition}\label{10Dec8}
\begin{enumerate}
\item The algebra $\mS_1$ is central.
 \item $Fa\neq 0$ and $aF\neq 0$ for all nonzero elements $a\in
 \mS_1$.
 \item $F$ is the smallest (with respect to inclusion) nonzero
 ideal of the algebra $\mS_1$ (i.e. $F$ is contained in all the
 nonzero ideals of $\mS_1$); $F^2=F$; $F$ is an essential left and
 right submodule of $\mS_1$; $F$ is the socle of the left and
 right $\mS_1$-module $\mS_1$.
 \item The set $\CJ (\mS_1)$ of all the ideals of the algebra
 $\mS_1$ is $\{ 0, F, \mS_1, F+a(K[x]+K[y])$ where $a=a(x)$ is a monic
 non-scalar polynomial of K[x] such that $a(0)\neq 0$\}; and two
 such ideals are equal,
 $F+a(K[x]+K[y])=F+b(K[x]+K[y])$,  iff $a=b$.
 \item $IJ=JI$ for all ideals $I$ and $J$ of the algebra $\mS_1$.
 \item $\spec (\mS_1)= \{ 0, F, \gm_a:=F+a(K[x]+K[y])$ where $a\in K[x]$ is a monic
 irreducible polynomial distinct from $x$\}. In particular,
 $\mS_1$ is a prime ring.
 \item $\Max (\mS_1)=\{ \gm_a\, | \, a$ is a monic
 irreducible polynomial distinct from $x$\}.
 \item Any proper ideal $I$ of the algebra $\mS_1$ such that
 $I\neq F$ is a unique finite  product of maximal ideals, i.e. $I=\prod
 \gm_a^{i_a}$ where all but finitely many natural numbers $i_a$ are equal to zero;  and $\prod \gm_a^{i_a}=\prod
 \gm_a^{j_a}$ iff all $i_a=j_a$ for all $a$.
 \item  $\prod \gm_a^{i_a}=\bigcap \gm_a^{i_a}$; $\prod \gm_a^{i_a}+\prod
 \gm_a^{j_a}=\prod \gm_a^{\min (i_a, j_a)}$ and $\prod \gm_a^{i_a}\bigcap \prod
 \gm_a^{j_a}=\prod \gm_a^{\max (i_a, j_a)}$. In particular, the
 lattice $\CJ (\mS_1)$ is distributive.
 \item The classical Krull dimension of the algebra $\mS_1$ is 2.
\end{enumerate}
\end{proposition}

{\it Remark}. For $K=\mathbb{C}$, statement 6 above and the fact
that $F$ is the minimal nonzero ideal of $\mS_1$ were proved by
Jordan (Corollary 7.6, \cite{Jordan-Down-up}) using a different
method.

$\noindent $

{\it Proof}. 1. By (\ref{Cenxy}), $Z(\mS_1)= \Cen (x) \cap \Cen
(y) = K$.

2. This follows at once from (\ref{mS1d}), (\ref{xyEij}), and
(\ref{xyEij1}).

3. $F^2=F$ since $F= M_\infty (K)$. The fact that $F$ is the
smallest nonzero ideal follows from statement 2. By statement 2,
$F$ is an essential left and right $\mS_1$-submodule of $\mS_1$.
Then $F$ is the socle of the module ${}_{\mS_1}\mS_1$ since the
$\mS_1$-module $F$ is semi-simple (see the proof of Lemma
\ref{c26Nov8}). Using the involution $\eta$ we see that $F$ is the
socle of the right $\mS_1$-module $\mS_1$.

4. Let $I$ be an ideal of the algebra $\mS_1$ which is distinct
from the ideals $0$, $F$ and $\mS_1$. Then $F\subset I$, and, by
(\ref{mS1d1}), $I=F+aK[x]$ for a monic non-scalar polynomial of
$K[x]$ with $a(0)\neq 0$ (since $I\neq \mS_1$). The rest is
obvious due to (\ref{mS1d1}).

5. This follows from statements 3 and  4 ($yx=1$, $xy= 1-E_{00}$,
$E_{00}\in F$).

6. Let $I$ and $J$ be nonzero ideals of the algebra $\mS_1$. Then
both of them contains the ideal $F$, by statement 3. Then
$IJ\supseteq F^2=F\neq 0$. Therefore, 0 is a prime ideal, i.e. the
algebra $\mS_1$ is  prime. $F$ is a prime ideal since
$\mS_1/F\simeq K[x,x^{-1}]$ is a Laurent polynomial algebra over a
field (a commutative domain).  From this fact it follows that the
ideal $F+a(K[x]+K[y])$ from statement 4 is prime iff in addition
the polynomial $a$ is irreducible. This observation finishes the
proof of statement 6.

7. Statement 7 follows from statement 6.

8. Statement 8 follows from statement 7.

9. Statement 9 follows from statement 8.

10. It follows from the inclusion of prime ideals $0\subset
F\subset \gm_a$ that $\clKdim (\mS_1)=2$ since the ideals $\gm_a$
are maximal. $\Box $

$\noindent $

By Proposition \ref{10Dec8}.(7), the map
$$ \Max(K[x,x^{-1}])\ra \Max (\mS_1), \;\; \gm \mapsto F+\mS_1\gm
= F+\gm \mS_1 , $$ is a bijection.

Let $n\geq2$. By (\ref{mS1d}), each element $a\in \mS_n= \mS_1\t
\mS_{n-1}$ has a {\em unique} presentation 
\begin{equation}\label{asum}
a=\l +\sum_{i\geq 1}(y^i\t \l_{-i}+x^i\t \l_i)+\sum_{i,j\geq
0}E_{ij}\t \l_{ij}
\end{equation}
for some elements $\l , \l_{\pm i}, \l_{ij}\in \mS_{n-1}$.  The
next lemma is crucial in the proofs of Theorems \ref{8Dec8} and
\ref{b23Dec8}).

\begin{lemma}\label{c23Dec8}
Let $I$ and $J$ be ideals of the algebra $\mS_n= \mS_1\t
\mS_{n-1}$, $n\geq 2$. Then
\begin{enumerate}
\item $I\cap (F\t \mS_{n-1})= F\t I_{n-1}$ for a unique ideal
$I_{n-1}$ of the algebra $\mS_{n-1}$. \item The ideal $I_{n-1}$ of
$\mS_{n-1}$ is the $K$-linear span in $\mS_{n-1}$ of the
coefficients $\l , \l_{\pm i}, \l_{ij}\in \mS_{n-1}$ in
 (\ref{asum}) for all the elements $a$ of the ideal $I$. \item If
$I\subseteq J$ then $I_{n-1}\subseteq J_{n-1}$. \item $(IJ)_{n-1}
= I_{n-1}J_{n-1}$.
\end{enumerate}
\end{lemma}

{\it Proof}. 1 and 2. The uniqueness in statement 1 is obvious (if
$F\t \ga = F\t \gb$ for two ideals $\ga$ and $\gb$ of $\mS_{n-1}$
then $\ga = \gb$, and vice versa). Let $I'_{n-1}$ be the
$K$-linear span in statement 2. Then $I'_{n-1}$ is an ideal of the
algebra $\mS_{n-1}$ and
$$ I\cap (F\t \mS_{n-1})\subseteq F\t I'_{n-1}.$$
To finish the proof of statements 1 and 2 it suffices to show that
 the reverse inclusion holds. Let $a\in I$, and so we have the
decomposition (\ref{asum}) for the element $a$. First, let us show
that $E_{kl}\t \l \in I$ for some natural numbers $k$ and $l$. Fix
sufficiently large natural numbers $k$ and $l$. Then
$E_{kl}(\sum_{i,j\geq 0} E_{ij}\t \l_{ij})E_{ll}=0$, and so
$$ E_{kl}aE_{ll} = E_{kl}\l +\sum_{i\geq 1} (E_{k, l+i}E_{ll}\t
\l_{-i}+E_{k,l-i}E_{ll}\t \l_i)= E_{kl}\t \l.$$ Similarly, for
sufficiently large natural numbers $k$ and $l$, and for all
natural numbers $i\geq 1$, $ E_{kl}y^iaE_{ll}= E_{kl}\t \l_i $ and
$E_{kl}a x^i E_{ll}= E_{kl}\t\l_{-i}$.  For all natural numbers
$i$, $E_{ii}aE_{ii}= E_{ii}\t (\l +\l_{ii})$ and so $E_{ii}\t
\l_{ii}\in I$ since $E_{ii}\t\l = E_{ik}(E_{kl}\t \l )E_{li}\in
I$. For all natural numbers $i>j$, $E_{ii}aE_{jj}= E_{ij}\t
(\l_{i-j} +\l_{ij})$ and $E_{jj}aE_{ii}= E_{ji}\t (\l_{-(i-j)}
+\l_{ji})$, and so $E_{ij}\t \l_{ij}, E_{ji}\t\l_{ji}\in I$. This
finishes the proof of statements 1 and 2.

 3. Statement 3 is obvious.

 4. Statement 4 follows from statements 1 and 2. By statement 2,
 $(IJ)_{n-1} \subseteq I_{n-1}J_{n-1}$. By statement 1, we have
 the inverse inclusion:
\begin{eqnarray*}
 F\t ( I_{n-1}J_{n-1}) &=& (F\t I_{n-1})(F\t J_{n-1})= I\cap (F\t
\mS_{n-1})\,  \cdot \, J\cap (F\t \mS_{n-1})\\
& \subseteq & (IJ)\cap (F\t \mS_{n-1}) = F\t (IJ)_{n-1},
\end{eqnarray*}
 and so $I_{n-1}J_{n-1}\subseteq (IJ)_{n-1}$.   $\Box $

$\noindent $

\begin{theorem}\label{8Dec8}
The set $\CJ (\mS_n)$ of ideals of the algebra $\mS_n$ satisfies
the ascending chain condition (the a.c.c., for short).
\end{theorem}

{\it Proof}. Recall that the set $\CJ (\mS_n)$ satisfies the
a.c.c. if each ascending chain of ideals in $\mS_n$ stabilizes,
i.e. has a largest element.  We use induction on $n$. The case
$n=1$ follows from the description of the set $\CJ (\mS_1)$
(Proposition \ref{10Dec8}.(4)). Suppose that $n\geq 2$ and that
the result is true for all $n'$ such that $n'<n$. Note that for
any algebra $A$, the set $\CJ (A)$ of its ideals satisfies the
a.c.c. iff the $A$-bimodule $A$ is Noetherian. Recall the
following (easy) generalization of the Hilbert Basis Theorem: {\em
An $A$-module $M$ is Noetherian iff the $K[x]\t A$-module $K[x]\t
M$ is Noetherian}. By induction, the $\mS_{n-1}$-bimodule
$\mS_{n-1}$ is Noetherian, hence the $K[x]\t\mS_{n-1}$-bimodule
$K[x]\t\mS_{n-1}$ is Noetherian, hence
$K[x,x^{-1}]\t\mS_{n-1}$-bimodule $K[x,x^{-1}]\t\mS_{n-1}$ is
Noetherian. Note that $F\t \mS_{n-1}$ is an ideal of the algebra
$\mS_n = \mS_1\t \mS_{n-1}$ such that $\mS_n /( F\t
\mS_{n-1})\simeq K[x,x^{-1}]\t\mS_{n-1}$ is a Noetherian
$K[x,x^{-1}]\t\mS_{n-1}$-bimodule, or, equivalently, a Noetherian
$\mS_n$-bimodule. For any ideal $I$ of $\mS_n$, by Lemma
\ref{c23Dec8}.(1),
$$ I\cap (F\t \mS_{n-1}) = F\t I_{n-1},$$
for some ideal $I_{n-1}$ of the algebra $\mS_{n-1}$.
 Therefore, the
$\mS_n$-bimodule $F\t \mS_{n-1}$ is Noetherian. It follows from
the short exact sequence of $\mS_n$-bimodules: 
\begin{equation}\label{FSI1}
0\ra F\t \mS_{n-1}\ra \mS_n = \mS_1\t \mS_{n-1} \ra K[x,x^{-1}]\t
\mS_{n-1} \ra 0
\end{equation}
that the $\mS_n$-bimodule $\mS_n$ is Noetherian since the
$\mS_n$-bimodules $ F\t \mS_{n-1}$ and $K[x,x^{-1}]\t \mS_{n-1}$
are Noetherian. This proves that the set $\CJ (\mS_n)$ satisfies
the a.c.c.. $\Box $

$\noindent $

{\it Definition}. For an algebra $A$ we say that its ideals {\em
commute} if $IJ = JI$ for all ideals $I$ and $J$ of the algebra
$A$.


\begin{theorem}\label{b23Dec8}
$IJ=JI$ for all ideals $I$ and $J$ of the algebra $\mS_n$.
\end{theorem}

{\it Proof}. To prove the result we use induction on $n$. The case
 $n=1$ is  Proposition \ref{10Dec8}.(5). So, let
$n>1$ and we assume that the result holds for all $n'<n$. By Lemma
\ref{c23Dec8}.(1), $I\cap (F\t \mS_{n-1}) = F\t I_{n-1}$ for some
ideal $I_{n-1}$ of the algebra $\mS_{n-1}$. Using (\ref{FSI1}), we
have the short exact sequence of $\mS_n$-modules
$$ 0\ra F\t I_{n-1} \ra I \ra \bI \ra 0$$
where $\bI$ is an ideal of the algebra $K[x,x^{-1}]\t \mS_{n-1}$
which is the image of the ideal $I$ under the epimorphism
$\mS_n\ra K[x,x^{-1}]\t \mS_{n-1}$. It is obvious that two ideals
$I$ and $I'$ of the algebra $\mS_n$ are equal iff $I_{n-1}=
I'_{n-1}$ and $\bI = \bI'$. Note that $(IJ)_{n-1} =
I_{n-1}J_{n-1}= J_{n-1} I_{n-1} = (JI)_{n-1}$ (by Lemma
\ref{c23Dec8}.(4)), and $\overline{IJ}= \bI\cdot  \bJ = \bJ\cdot
\bI = \overline{JI}$ (by the same arguments and induction).
Therefore, $IJ = JI$. $\Box $

$\noindent $

{\bf The associated graded algebra ${\rm gr} (\mS_n)$ and the
algebra $\mrD_n$}.

$\noindent $

{\it Definition}. The 
algebra $\mathrm{D}_n$ is an algebra generated over a field $K$
by $2n$ elements $x_1, \ldots , x_n, y_n, \ldots , y_n$ that
satisfy the defining relations:
$$ y_1x_1=\cdots = y_nx_n=0 , \;\; [x_i, y_j]=[x_i, x_j]= [y_i,y_j]=0
\;\; {\rm for\; all}\; i\neq j.$$ Clearly, $\mathrm{D}_n\simeq
\mathrm{D}_1^{\t n}$ and $\mrD_n =\bigoplus_{\alpha, \beta \in
\N^n} Kx^\alpha y^\beta$. The canonical generators $x_i$, $y_j$
$(1\leq i,j\leq n)$ of the algebra $\mrD_n$ determine the
ascending filtration $\{ \mrD_{n, \leq i}\}_{i\in \N}$ on the
algebra $\mrD_n$ where $\mrD_{n, \leq i}:= \bigoplus_{|\alpha
|+|\beta |\leq i} Kx^\alpha y^\beta$. Then  $\dim (\mrD_{n,\leq
i})={i+2n \choose 2n}$ for $i\geq 0$, and so $\GK (\mrD_n )=2n$.

The associated graded algebra ${\rm gr} (\mS_n) := \bigoplus_{i\in
\N} \mS_{n, \leq i}/ \mS_{n, \leq i-1}$ is isomorphic to the
algebra $\mathrm{D}_n$. The {\em nil-radical} $\gn = \gn
(\mathrm{D}_n)$ of the algebra $\mathrm{D}_n$ (i.e. the sum of all
the nilpotent ideals of $\mathrm{D}_n$) is equal to
$\bigoplus_{\exists i : \alpha_i\beta_i\neq 0}Kx^\alpha y^\beta$,
and $\gn^{n+1}=0$. The factor algebra $\mathrm{D}_n/\gn\simeq
\bigotimes_{i=1}^n K[\bx_i, \by_i]/(\by_i\bx_i)$ is the tensor
product of the commutative algebras $ K[\bx_i,
\by_i]/(\by_i\bx_i)$.

The involution $\eta$ of the algebra $\mS_n$ respects the
filtration $\{ \mS_{n, \leq i}\}_{i\in \N}$, i.e. $\eta (\mS_{n,
\leq i})=\mS_{n, \leq i}$ for all $i\geq 0$; and so the associated
graded algebra ${\rm gr} (\mS_n)$ inherits the involution
$$\eta :\mathrm{D}_n\ra
\mathrm{D}_n, \;\; x_i\mapsto y_i, \;\; y_i\mapsto x_i, \;\;
i=1,\ldots , n. $$ In particular, the algebra $\mathrm{D}_n$ is
self-dual. The algebra $\mathrm{D}_n$ is neither left nor right
Noetherian as it contains the infinite direct sum
$\bigoplus_{\beta \in \N^n} \mathrm{D}_nx_1\cdots x_ny^\beta$ of
nonzero left ideals. The simple $\mathrm{D}_n$-modules can be
easily described,
$$\widehat{\mathrm{D}_n}=\widehat{\mathrm{D}_n/\gn }.$$ In particular, all the simple  $\mathrm{D}_n$-modules
are finite dimensional. The prime and maximal spectra of the
algebra $\mathrm{D}_n$ are easily found since $\Spec
(\mathrm{D}_n) = \Spec(\mathrm{D}_n/\gn )$ and $\Max
(\mathrm{D}_n) = \Max (\mathrm{D}_n/ \gn )$.


\section{Classification of simple $\mS_n$-modules}

In this section, we classify all the simple $\mS_n$-modules
(Theorem \ref{17Dec8}.(1)). It is proved that for each simple
$\mS_n$-module $M$, its endomorphism algebra $\End_{\mS_n}(M)$ is
a finite {\em field} extension of the field $K$ (Theorem
\ref{17Dec8}.(4)), and the multiplicity $e(M)$ of $M$ is equal to
$\dim (\End_{\mS_n}(M))$ (Corollary \ref{d17Dec8}). This is the
second instance known to me when the multiplicity of a simple
module is equal to the dimension of its endomorphism division
algebra. In  \cite{Bav-holmodp} this was proved for certain simple
modules over the ring $\CD (P_n)$ of differential operators on the
polynomial algebra $P_n$ over a perfect field of {\em prime}
characteristic. Note that the algebra $\CD (P_n)$ is neither left
nor right Noetherian and   not finitely generated either.

$\noindent $

{\bf The algebra $\mS_n$ is $\Z^n$-graded}. The algebra $\mS_n
=\bigoplus_{\alpha \in \Z_n}\mS_{n,\alpha}$ is a $\Z^n$-graded
algebra where $\mS_{n,\alpha}:= \mS_{1,\alpha_1}\t\cdots \t
\mS_{1,\alpha_n}$, $\alpha = (\alpha_1, \ldots , \alpha_n)$,
$$ \mS_{1,i}:=
\begin{cases}
x^i\mS_{1,0}= \mS_{1,0}x^i & \text{if $i\geq 1$ },\\
\mS_{1,0}& \text{if $i=0$ },\\
y^{|i|}\mS_{1,0}= \mS_{1,0}y^{|i|} & \text{if $i\leq -1$ ,}
\end{cases}
$$
$\mS_{1,0}:= K\langle E_{00}, E_{11}, \ldots \rangle= K\oplus
KE_{00}\oplus KE_{11}\oplus\cdots$ is a commutative non-Noetherian
algebra $( KE_{00} \subset  KE_{00}\oplus KE_{11}\subset \cdots$
is an ascending chain of ideals of the algebra $\mS_{1,0}$).

The polynomial algebra $P_n=\bigoplus_{\alpha \in \N^n} P_{n,
\alpha}$ is an $\N^n$-graded algebra (and, automatically,
 a $\Z^n$-graded algebra) where $P_{n,\alpha}:=Kx^\alpha$. Moreover,
$\mS_{n,\alpha}P_{n,\beta}\subseteq P_{n, \alpha +\beta}$ for all
$\alpha , \beta \in \Z^n$. Therefore, the polynomial algebra $P_n$
is a $\Z^n$-graded $\mS_n$-module.

Note that the involution $\eta$ reverses the $\Z^n$-grading of the
algebra $\mS_n$, i.e.
$$ \eta (\mS_{n, \alpha}) = \mS_{n,-\alpha }, \;\; \alpha \in
\Z^n,$$ and it acts as the identity map on the algebra
$\mS_{n,0}$.

For an algebra $A$, let $\hA$ denote the set of isoclasses of
simple $A$-modules. For a simple $A$-module $M$, $[M]$ is its
isoclass. We usually drop the brackets $[,]$ if this does not lead
to confusion.

$\noindent $

{\bf The simple $\mS_1$-modules}. The next Lemma gives all the
simple $\mS_1$-modules.

\begin{lemma}\label{a17Dec8}
\begin{enumerate}
\item $\hmS_1= \{ K[x], \mS_1/ \gm_a= \mS_1/(F+a(K[x]+K[y]))$
where $\gm_a\in \Max (\mS_1)\}$ and $\hmS_1 \simeq \{ K[x]\} \cup
\widehat{K[x,x^{-1}]}$.
 \item $\End_{\mS_1}(K[x])\simeq K$ and $\End_{\mS_1}(\mS_1/
 \gm_a)\simeq K[x,x^{-1}]/(a)$ for all $\gm_a\in \Max (\mS_1)$.
\end{enumerate}
\end{lemma}

{\it Proof}. 1. Let $M$ be a simple $\mS_1$-module. Then either
$FM=M$ or $FM=0$. In the first case, $E_{ij}m\neq 0$ for some
matrix unit $E_{ij}$ and a nonzero element $m\in M$. Recall that
${}_{\mS_1}\mS_1E_{ij} = \mS_1E_{0j}\simeq K[x]$ (a simple
$\mS_1$-module). Since the $\mS_1$-homomorphism $\mS_1E_{ij}\ra
\mS_1E_{ij}m$, $a\mapsto am$, is nonzero, $M=\mS_1E_{ij}m\simeq
\mS_1E_{ij}\simeq K[x]$.

In the second case, $FM=0$, the $\mS_1$-module $M$ is in fact the
module over the factor algebra $\mS_1/ F\simeq K[x,x^{-1}]$, hence
$M\simeq \mS_1/ \gm_a$ for some maximal ideal $\gm_a\in \Max
(\mS_1)$ (see Proposition \ref{10Dec8}.(7)). It is obvious that
$\hmS_1 \simeq \{ K[x]\} \cup \widehat{K[x,x^{-1}]}$.

2. Since ${}_{\mS_1}K[x]\simeq \mS_1/ \mS_1y$ and $\ann_{K[x]}(y)
= K$, we have $\End_{\mS_1}(K[x])\simeq K$. It is obvious that
 $\End_{\mS_1}(\mS_1/ \gm_a)\simeq \End_{K[x,x^{-1}]}(K[x,x^{-1}]/ (a))\simeq
 K[x,x^{-1}]/(a)$. $\Box $

$\noindent $

{\bf The simple $\mS_n$-modules}. Our next goal is to generalize
Lemma \ref{a17Dec8} for an arbitrary $n$. First, we introduce more
notation. For a subset $\CN=\{ i_1, \ldots , i_s\}$ of the set of
indices $\{ 1, \ldots , n\}$, let $C\CN$ be its complement, $|\CN
|=s$, $\mS_\CN := \mS_1(i_1) \t\cdots \t \mS_1(i_s)$,
\begin{equation}\label{gmCN}
\ga_\CN :=F\t\mS_1(i_2) \t\cdots \t \mS_1(i_s)+\cdots +\mS_1(i_1)
\t\cdots \t \mS_1(i_{s-1})\t F,
\end{equation}
 $P_\CN := K[x_{i_1}, \ldots , x_{i_s}]$.
Clearly, $\mS_n= \mS_\CN \t \mS_{C\CN }$. Let $L_\CN:= K[x_{i_1},
x_{i_1}^{-1}, \ldots , x_{i_s}, x_{i_s}^{-1}]$. Then $\mS_\CN /
\ga_\CN \simeq L_\CN$. The $\mS_\CN$-module $P_\CN$ is simple with
$\End_{\mS_\CN }(P_\CN )\simeq K$. For each maximal ideal $\gm $
of the algebra $L_{C\CN}$, the $\mS_{C\CN}$-module $L_{C\CN} / \gm
\simeq \mS_{C\CN} / (\ga_{C\CN} +\mS_{C\CN} \gm )$ is simple.  Let
$F_{C\CN}=\bigotimes_{i\in C\CN }F(i) = F^{\t (n-s)}$. Hence, the
tensor product 
\begin{equation}\label{NCNm}
M_{\CN, \gm}:= P_\CN \t (L_{C\CN}/\gm )
\end{equation}
is a simple $\mS_n$-module. The annihilator $\ann_{\mS_n}(M_{\CN,
\gm})$ is equal to $\mS_\CN\t (F_{C\CN}+\mS_{C\CN}\gm )$.
Therefore, two such modules are isomorphic, $M_{\CN_1,
\gm_1}\simeq M_{\CN_2, \gm_2}$, iff $\CN_1=\CN_2$ and
$\gm_1=\gm_2$ iff $\ann_{\mS_n}(M_{\CN_1,
\gm_1})=\ann_{\mS_n}(M_{\CN_2, \gm_2})$. Let 
\begin{equation}\label{hSnN}
\hmS_{n,\CN } := \{ M_{\CN, \gm}\, | \, \gm \in \Max (L_{C\CN}
)\}.
\end{equation}
In particular, $\hmS_{n, \{ 1, \ldots , n\} }=\{ P_n\}$ and
$\hmS_{n, \emptyset }= \{ L_n/ \gm \; | \; \gm \in \Max (L_n)\}$.
The subsets $\{ \hmS_{n, \CN } \}$ of $\hmS_n$ are disjoint.

\begin{theorem}\label{17Dec8}
\begin{enumerate}
\item $\hmS_n = \coprod_{\CN \subseteq \{ 1, \ldots , n\}}\hmS_{n,
\CN }$, a disjoint union.  \item $\End_{\mS_n}(P_n) \simeq K$,
$\ann_{\mS_n}(P_n)=0$, and $\GK (P_n) = n$.  \item
$\End_{\mS_n}(M_{\CN , \gm}) \simeq L_{C\CN} /\gm $,
$\ann_{\mS_n}(M_{\CN , \gm})=\mS_\CN \t (F_{C\CN} +\mS_{C\CN} \gm
)$, and $\GK (M_{\CN , \gm}) = |\CN |$. \item The endomorphism
algebra of each simple $\mS_n$-module is a finite field extension
of $K$.
\end{enumerate}
\end{theorem}

{\it Proof}. 1. We use induction on $n$. The case $n=1$ is Lemma
\ref{a17Dec8}. Suppose that $n>1$ and that the result is true for
all $n'<n$. Let $M$ be a simple $\mS_n$-module. Then either $\ga_n
M=0$ or $\ga_nM=M$. In the first case, the module $M$ is a simple
$(\mS_n/ \ga_n=L_n)$-module, and so $M\in \hmS_{n, \emptyset}$.

In the second case, $F(i) M\neq 0$ for some $i\in \{ 1, \ldots ,
n\}$ where $F(i)$ is the smallest nonzero ideal of the algebra
$\mS_1(i)$. Without loss of generality we may assume that $i=n$.
Then $\mS_1(n)$-module $M$ contains a simple $\mS_1(n)$-submodule
isomorphic to the $\mS_1(n)$-module $K[x_n]$. Since $\mS_n =
\mS_{n-1} \t \mS_1(n)$ and $\End_{\mS_1(n)}(K[x_n])\simeq K$, we
have $M\simeq N\t K[x_n]$ for a simple $\mS_{n-1}$-module $N$.
Now, induction completes the argument.

2 and 3. Statements 2 and 3 are obvious.

4. Statement 4 follows from statements 2 and 3.  $\Box $

$\noindent $

An algebra is called a {\em primitive} algebra provided there
exists a simple faithful module. The {\em Jacobson radical} of an
algebra is the largest ideal that annihilates all the simple
modules.
\begin{corollary}\label{a19Dec8}
The polynomial algebra $P_n$ is the only faithful simple
$\mS_n$-module (and so the Jacobson radical of $\mS_n$ is zero,
and the algebra $\mS_n$ is primitive).
\end{corollary}

{\it Proof}. The polynomial algebra $P_n$ is a faithful simple
$\mS_n$-module. The fact that it is the only one follows from
Theorem \ref{17Dec8}.(3). $\Box $

$\noindent $

Note that the annihilator of simple module is a prime ideal. The
following corollary shows that the annihilators separate the
simple modules.

\begin{corollary}\label{b17Dec8}
The map $\hmS_n\ra \spec (\mS_n)$, $M\mapsto \ann_{\mS_n}(M)$, is
an injection.
\end{corollary}

{\it Proof}. We have seen above that simple modules are isomorphic
iff their annihilators coincide. $\Box $

$\noindent $

For the Weyl algebra $A_n$ over a field $K$ of characteristic
zero, the {\em Inequality of Bernstein} says that $\GK (M)\geq n$
for all nonzero finitely generated $A_n$-modules $M$. Note that
$\GK (A_n) = \GK (\mS_n) = 2n$. The corollary  below shows that
for simple $\mS_n$-modules the `opposite' Inequality of Bernstein
is true.

\begin{corollary}\label{c17Dec8}
$\{ \GK (M) \; | \; M\in \hmS_n\}= \{ 0,1, \ldots , n\}$.
\end{corollary}

{\it Proof}. Theorem \ref{17Dec8}.(2,3). $\Box $

$\noindent $

Recall that the algebras $\mS_n$ and $P_n$ are equipped with the
standard filtrations $\{ \mS_{n, \leq i}\}_{i\in \N}$ and $\{
P_{n, \leq i}\}_{i\in \N}$. For each simple $\mS_n$-module $
M_{\CN ,  \gm}$ (see (\ref{NCNm})), let $1\t \overline{1}$ be its
canonical generator where $1\in P_\CN$ and $ \overline{1}\in
L_{C\CN }/\gm$. The $\mS_n$-module $M_{\CN ,  \gm}$ admits the
standard filtration
$$ \{ M_{\CN,  \gm , \leq i}:= \mS_{n, \leq i} \, 1\t
\overline{1}\}.$$ There exists a natural number, say $j$, such
that
$$  P_{\CN, \leq i-j}\t (L_{C\CN } / \gm ) \subseteq M_{\CN,  \gm , \leq i}
\subseteq P_{\CN, \leq i}\t (L_{C\CN } / \gm ), \;\; i\gg 0,$$
 and so the Hilbert function of the module $M_{\CN
, \gm}$ has the form:
$$ \dim (M_{\CN ,  \gm , \leq i})= \frac{\dim(L_{C\CN } /\gm )}{|\CN |!}\,
i^{|\CN |}+\cdots , \;\; i\geq 0, $$ where the three dots denote
smaller terms. It follows that the {\em multiplicity} $e(M_{\CN ,
\gm})$ of the $\mS_n$-module $M_{\CN , \gm}$ exists and  is equal
to
$$ \dim (L_{C\CN } / \gm ) = \dim (\End_{\mS_n}(M_{\CN
,  \gm})),$$ by Theorem \ref{17Dec8}.(2,3).

\begin{corollary}\label{d17Dec8}
For each simple $\mS_n$-module $M$, its multiplicity $e(M)$ is
equal to $\dim (\End_{\mS_n}(M))$.
\end{corollary}


\section{The prime and maximal spectra of the algebra $\mS_n$}

In this section, the prime and  maximal spectra of the algebra
$\mS_n$ are found (Theorems \ref{21Dec8} and \ref{d21Dec8}); it is
proved that the classical Krull dimension of the algebra $\mS_n$
is $2n$ (Theorem \ref{f21Dec8}), and the algebra $\mS_n$ is a
central, prime, catenary algebra (Proposition \ref{b21Dec8},
Theorem \ref{23Dec8}). Formulae for the height and the co-height
of primes of the algebra $\mS_n$ are found via combinatorial data
(Theorem \ref{22Dec8}). In many arguments, the height 1 primes of
the algebra $\mS_n$ play an important role. Their classification
is given in Lemma \ref{a21Dec8}.

$\noindent $

{\bf The algebra $\mS_n$ is a prime algebra}. Recall that $F_n:=
F^{\t n }$ is an ideal of the algebra $\mS_n$. It follows from
$F_n=\bigoplus_{\alpha , \beta \in \N^n}KE_{\alpha \beta}$ that
$F_n$ is a {\em simple} $\mS_n$-bimodule.

\begin{proposition}\label{b21Dec8}
\begin{enumerate}
\item The algebra $\mS_n$ is central. \item  $F_n a \neq 0$ and
$aF_n\neq 0$ for all nonzero elements $a$ of the algebra $\mS_n$.
\item $F_n$ is the smallest (with respect to inclusion) nonzero
ideal of the algebra $\mS_n$ (i.e. $F_n$ is contained in all
nonzero ideals of $\mS_n$); $F_n^2= F_n$; $F_n$ is an essential
left and right submodule of $\mS_n$; $F_n$ is the socle of the
left and right $\mS_n$-module $\mS_n$; $F_n$ is the socle of the
$\mS_n$-bimodule $\mS_n$ and $F_n$ is a simple $\mS_n$-bimodule.
\item The algebra $\mS_n$ is a prime algebra.\item Every nonzero
ideal of the algebra $\mS_n$ is an essential left and right
submodule of $\mS_n$.
\end{enumerate}
\end{proposition}

{\it Proof}. 1. The fact that the algebra $\mS_n$ is central
follows from (\ref{Cenxy}) and Proposition \ref{10Dec8}.(1).

2. For $n=1$, this is Proposition \ref{10Dec8}.(2). Suppose that
$n>1$. Let $a$ be a nonzero element of the algebra $\mS_n =
\mS_1\t \mS_{n-1}$. The element  $a$ is a sum $\sum a_i\t b_i$ for
some linearly independent elements $a_i \in \mS_1$ and some
linearly independent elements $b_i\in \mS_{n-1}$. Now, the result
follows from Proposition \ref{10Dec8}.(2) since $Fa_1\neq 0 $ and
$a_1F\neq 0$.

3. The fact that $F_n$ is the smallest nonzero ideal of the
algebra $\mS_n$ follows from statement 2 and the fact that $F_n$
is a simple $\mS_n$-bimodule. $F_n^2= (\bigoplus_{\alpha , \beta
\in \N^n} KE_{\alpha \beta})^2=\bigoplus_{\alpha , \beta \in \N^n}
KE_{\alpha \beta}= F_n$. By statement 2, $F_n$ is an essential
left and right submodule of $\mS_n$. Note that
$$ {}_{\mS_n} F_n = \bigoplus_{\alpha \in \N^n} (\mS_n E_{0,
\alpha_1}\cdots E_{0, \alpha_n}) \simeq \bigoplus_{\alpha \in
\N^n} P_n,$$ and so the left $\mS_n$-module $F_n$ is semi-simple.
Therefore, $F_n$ must be the socle of the left $\mS_n$-module
$\mS_n$ since $F_n$ is an essential submodule of $\mS_n$. Applying
the involution $\eta$ and using the fact that $\eta (F_n) = F_n$,
we see that $F_n$ is the socle of the right $\mS_n$-module
$\mS_n$.

4. The algebra $\mS_n$ is a prime algebra since $F_n$ is the
smallest  nonzero ideal of the algebra $\mS_n$ and $F_n^2= F_n$
(let $I$ and $J$ be nonzero ideals of $\mS_n$; then $IJ\supseteq
F_n^2= F_n\neq 0$).

5. Statement 5 follows from the fact that $F_n$ is the smallest
nonzero ideal of $\mS_n$ and $F_n$ is an essential left and right
submodule of $\mS_n$ (statement 2). $\Box $

$\noindent $

{\bf The set of height 1 primes of $\mS_n$}.   Consider the ideals
of the algebra $\mS_n$:
$$\gp_1:=F\t \mA_{n-1},\; \gp_2:= \mS_1\t F\t \mS_{n-2}, \ldots ,
 \gp_n:= \mS_{n-1} \t F.$$ Then $\mS_n/\gp_i\simeq
\mS_{n-1}\t (\mS_1/F) \simeq  \mS_{n-1}\t K[x_i, x_i^{-1}]$ and
$\bigcap_{i=1}^n \gp_i = \prod_{i=1}^n \gp_i =F^{\t n }$. Clearly,
$\gp_i\not\subseteq \gp_j$ for all $i\neq j$.

Let $\gp$ be a prime ideal of an algebra $A$. Let $\hht (\gp
)=\hht_A(\gp)$ denote the {\em height} of the prime ideal $\gp$
and  $\cht (\gp ) = \cht_A(\gp ) := \clKdim (A/ \gp )$ denote  the
{\em co-height} of the prime ideal $\gp$.

The next lemma  gives all the height 1 primes of the algebra
$\mS_n$. Surprisingly,  there are only finitely many of them
(bearing in mind that $\mS_n$ is a prime algebra of classical
Krull dimension $2n$ (Theorem \ref{f21Dec8}) and $\mS_n$ is a
`left localization' of the polynomial algebra $P_n$).

\begin{lemma}\label{a21Dec8}
The set of height 1 prime ideals of the algebra $\mS_n$ is $\{
\gp_1, \ldots , \gp_n\}$.
\end{lemma}

{\it Proof}. The algebras $\mS_{n-1}$  and $K[x_i, x_i^{-1}]$ are
prime algebras, then so is their tensor product $\mS_{n-1}\t
K[x_i, x_i^{-1}]\simeq \mS_n / \gp_i$ (Proposition
\ref{c21Dec8}.(1)), and so $\gp_i$ is a prime ideal. Clearly,
$\hht (\gp_i)\geq 1$ for all $i$. Suppose that  $\hht (\gp_i)> 1$
for some $i$, we seek a contradiction. Then there is a nonzero
prime $\gp$ which is strictly contained in the ideal $\gp_i$. Then
$\gp_i\supset \gp \supseteq F_n =\prod_{j=1}^n \gp_j$, and so
$\gp_i\supseteq \gp_j$ for some $j\neq i$, a contradiction.
Therefore, $\hht (\gp_i) =1$ for all $i$.  $\Box $

$\noindent $

The primes $\gp_1, \ldots , \gp_n$ play a prominent role in many
proofs of this paper. To find the prime spectrum $\spec (\mS_n)$
of the algebra $\mS_n$ we need the following result.

\begin{proposition}\label{c21Dec8}
Let $A$ be an algebra over the field $K$.
\begin{enumerate}
\item Then the algebra $\mS_n\t A$ is a prime algebra iff the
algebra $A$ is so. \item If $\gp$ is  a prime ideal of the algebra
$A$ then the ideal $\mS_n\t \gp$ of the algebra $\mS_n\t A$ is a
prime ideal, and vice versa.
\end{enumerate}
\end{proposition}

{\it Proof}. 1. Suppose that the algebra $A$ is not prime. Then
$\ga \gb =0$ for some nonzero ideals $\ga $ and $\gb$ of the
algebra $A$. The ideals $\mS_n\t \ga$ and $\mS_n\t \gb$ of the
algebra $\mS_n\t A$ are nonzero. Since their product is zero the
algebra $\mS_n\t A$ is not prime.

In order to finish the proof it suffices to show that if the
algebra $A$ is prime then so is the algebra $\mS_n\t A$. Let $\ga$
be a nonzero ideal of the algebra $\mS_n\t A$. Then $F_n\ga\neq
0$, by Proposition \ref{b21Dec8}.(2). Note that $F_n\ga \subseteq
\ga$. Let $u=E_{\alpha \beta} \t a+\cdots + E_{\s \rho}\t a'$ be a
nonzero element of $F_n\ga$ where $E_{\alpha\beta}, \ldots ,
E_{\s\rho}$ are distinct matrix units;  $a, \ldots , a'\in A$, and
$a\neq 0$. Then $0\neq E_{\alpha \beta} \t a= E_{\alpha\alpha}
uE_{\beta\beta} \in \ga$, and so $F_n\t AaA\subseteq \ga$.
Similarly, $F_n\t AbA\subseteq \gb$ for some nonzero element $b\in
 A$. Then
$$\ga\gb \supseteq F_n\t AaA\, \cdot \, F_n\t AbA= F_n\t (
AaA\cdot AbA)\neq 0$$ since $F_n=\bigoplus_{\alpha, \beta \in
\N^n}KE_{\alpha\beta}$ and $AaA\cdot AbA\neq 0$ ($A$ is a prime
algebra). Therefore, $\mS_n\t A$ is a prime algebra.

2. Statement 2 follows from statement 1 since $\mS_n\t A/(\mS_n\t
\gp ) \simeq \mS_n\t (A/ \gp )$. $\Box $

$\noindent $

 Let
$\ga_n:= \gp_1+\cdots +\gp_n$. Then the factor algebra
\begin{equation}\label{SnSn}
\mS_n/ \ga_n\simeq (\mS_1/F)^{\t n } \simeq \bigotimes_{i=1}^n
K[x_i, x_i^{-1}]= K[x_1, x_1^{-1}, \ldots , x_n, x_n^{-1}]=:L_n
\end{equation}
is a  Laurent polynomial algebra in $n$ variables,  and so $\ga_n$
is a prime ideal of co-height $n$ of the algebra $\mS_n$. The
algebra $L_n$ is commutative, and so $$ [a,b]\in \ga_n\;\; {\rm
for \; all}\;\; a,b\in \mS_n .$$ Since $\eta (\ga_n) = \ga_n$, the
involution of the algebra $\mS_n$ induces the {\em automorphism}
$\overline{\eta}$ of the factor algebra $\mS_n / \ga_n$ by the
rule:
$$\overline{\eta}: L_n\ra L_n, \;\; x_i\mapsto x_i^{-1}, \;\; i=1,
\ldots , n.$$ It follows that $\eta (ab) - \eta (a) \eta (b)\in
\ga_n$ for all elements $a, b\in \mS_n$. For each subset $\CN$ of
the  set $\{ 1, \ldots , n\}$, consider the epimorphism
\begin{equation}\label{pCN}
\pi_\CN : \mS_\CN \ra \mS_\CN/ \ga_\CN\simeq L_\CN, \;\; a\mapsto
a+ \ga_\CN,
\end{equation}
where $\ga_\CN$ is defined in (\ref{gmCN}). By Proposition
\ref{c21Dec8}.(2), there is the injection
$$ \spec (L_{C\CN}) \ra \spec (\mS_n) , \;\; \gq \mapsto \mS_\CN \t
\pi_{C\CN}^{-1}(\gq ).$$ The image of this injection is denoted by
$$ \spec (\mS_n, \CN ) :=\{ \mS_\CN \t
\pi_{C\CN}^{-1}(\gq ) \, | \, \gq \in \spec (L_{C\CN} )\}.$$
Clearly, $\spec (\mS_n, \CN )\cap \spec (\mS_n, \CM )=\emptyset$
if $\CN \neq \CM$ (see also Lemma \ref{e21Dec8} for details). Note
that $\spec (\mS_n , \emptyset )= \{ \pi^{-1}_{\{ 1, \ldots , n\}
} (\gq )\, | \, \gq \in \spec (L_n)\}\simeq \spec (L_n)$ and
$\spec (\mS_n , \{ 1, \ldots , n\} )= \{ 0\}$ since $\pi_\emptyset
:K\ra K$, $\l \mapsto \l$.

The next theorem shows that all the prime ideals of the algebra
$\mS_n$ can be obtained in this way.
\begin{theorem}\label{21Dec8}
\begin{enumerate}
\item $\spec (\mS_n) = \coprod_{\CN \subseteq \{ 1, \ldots , n\} }
\spec (\mS_n , \CN )$, the disjoint union.  \item Each  prime
ideal $\gp$ of the algebra $\mS_n$ can be uniquely written as
$\mS_\CN \t \pi^{-1}_{C\CN}(\gq )$ for some subset $\CN$ of the
set $\{ 1, \ldots , n\}$ and some prime ideal $\gq$ of the algebra
$L_{C\CN }$.
\end{enumerate}
\end{theorem}

{\it Proof}.  To prove the theorem it suffices to show that each
nonzero  prime ideal $\gp$ belongs to the union (statement 2
follows from statement 1, and uniqueness is obvious). To prove
this we use induction on $n$. The case $n=1$ is obvious
(Proposition \ref{10Dec8}.(6)). So, let $n>1$. By Lemma
\ref{a21Dec8}, $\gp_i\subseteq \gp$ for some $i$ since $\gp$ is a
nonzero prime ideal. Up to permutation of the indices $1, \ldots ,
n$, we may assume that $i=n$, i.e. $\gp_n\subseteq \gp$.
Therefore, $\overline{\gp } := \gp / \gp_n$ is a prime ideal of
the factor algebra $\mS_n / \gp_n\simeq \mS_{n-1} \t K[x_n,
x_n^{-1}]$. The elements of the algebra $K[x_n, x_n^{-1}]$ are
central in the algebra $\mS_{n-1} \t K[x_n, x_n^{-1}]$. The
algebra $\mS_{n-1} \t K[x_n, x_n^{-1}]$ is a subalgebra of the
algebra $\mS_{n-1} \t K(x_n)$ which is the localization of the
algebra $\mS_{n-1} \t K[x_n, x_n^{-1}]$ at the {\em central}
multiplicative set $S:=  K[x_n, x_n^{-1}]\backslash  \{ 0 \}$ all
the elements of which are regular. Note that the algebra
$\mS_{n-1} \t K(x_n)$ is the 
algebra $\mS_{n-1} (K(x_n))$
over the field $K(x_n)$ of rational functions. The localization
$S^{-1} \overline{\gp}= \overline{\gp }\t_{K[x_n,x_n^{-1}]}K(x_n)$
of the nonzero prime ideal $\overline{\gp}$ of the algebra
$\mS_{n-1} \t K[x_n, x_n^{-1}]$ is a nonzero prime ideal of the
algebra $\mS_{n-1} \t K(x_n)$ since the set $S$ is central. There
are two cases to consider: either $S^{-1} \overline{\gp}=
\mS_{n-1} \t K(x_n)$ or $S^{-1} \overline{\gp}\neq \mS_{n-1} \t
K(x_n)$.

In the first case, the ideal $\gp$ contains the prime ideal
$\mS_{n-1}\t \gm$ for some maximal ideal $\gm$ of the algebra
$\mS_1 (n)$ (see Proposition \ref{10Dec8}.(6,7)), and so $\gp /
(\mS_{n-1} \t \gm )$ is a prime ideal of the factor algebra
$$ \mS_n / (\mS_{n-1} \t\gm )\simeq (\mS_{n-1} \t  \mS_1 (n))/(\mS_{n-1}\t  \gm
)\simeq\mS_{n-1} \t (\mS_1(n)/ \gm ) \simeq \mS_{n-1} \t K_\gm$$
over the field $K_\gm := \mS_1 (n) / \gm$ which is a finite field
extension of the field $K$ (Proposition \ref{10Dec8}.(6)). The
algebra $\mS_{n-1} \t K_\gm$ is the 
algebra $\mS_{n-1}$ over the field $K_\gm$. By induction, $\gp /
(\mS_{n-1} \t \gm )= \mS_\CN \t \pi^{-1}_{C'\CN }(\gq )$ for some
subset $\CN$ of $\{ 1, \ldots , n-1\}$ and some  prime ideal $\gq$
of the algebra $L_{C'\CN}\t K_\gm$ where $C'\CN : = \{ 1, \ldots ,
n-1\} \backslash \CN$ and
$$ \pi_{C'\CN}:\mS_{C'\CN}\t K_\gm \ra L_{C'\CN}\t K_\gm. $$
Consider the commutative diagram of the algebras
$$\xymatrix{\mS_{C\CN}=\mS_{C'\CN}\t \mS_1(n)\ar[d]^f\ar[r]^{\pi_{C\CN}} & L_{C\CN}=L_{C'\CN}\t K[x_n,x_n^{-1}]\ar[d]^g\\
\mS_{C'\CN}\t K_\gm  \ar[r]^{\pi_{C'\CN}} &L_{C'\CN}\t K_\gm}$$
where all four maps are obvious epimorphisms (and $C\CN = \{ 1,
\ldots , n\} \backslash \CN$). Then $\gq':= g^{-1} (\gq )$ is a
prime ideal of the algebra $L_{C\CN }$. It follows from $ (1\t f)
(\gp ) = \gp / (\mS_{n-1} \t \gm ) = \mS_{\CN} \t
\pi^{-1}_{C'\CN}(\gq )$ (where $1\t f : \mS_n = \mS_\CN \t
\mS_{C\CN }\ra \mS_\CN \t \mS_{C'\CN } \t K_\gm $) that
$$ \gp = \mS_\CN \t f^{-1} \pi^{-1}_{C'\CN } (\gq ) = \mS_\CN \t
\pi^{-1}_{C\CN } (g^{-1}(\gq )) = \mS_\CN \t \pi^{-1}_{C\CN}
(\gq').$$ So, the result is true in the first case.

It remains to consider the second case when $S^{-1} \overline{\gp
}\neq \mS_{n-1} \t K(x_n)$. Note that $\mS_{n-1} \t K(x_n)$ is the
 algebra $\mS_{n-1} ( K(x_n))$ over the field $K(x_n)$. By
induction, $S^{-1} \overline{\gp}= \mS_\CN \t \pi^{-1}_{C'\CN}(\gq
)$ for some subset $\CN$ of the set $\{ 1, \ldots , n-1\}$ and
some prime ideal $\gq$ of the algebra $L_{C'\CN}\t K(x_n)$ where
$C'\CN = \{ 1, \ldots , n-1\}\backslash \CN$ and $$
\pi_{C'\CN}:\mS_{C'\CN}\t K(x_n) \ra L_{C'\CN}\t K(x_n). $$
Consider the commutative diagram of the algebras
$$\xymatrix{\mS_{n-1}\t K[x_n, x_n^{-1}]\ar[d]^u\ar[r]^{\v} & \mS_{n-1}\t K(x_n)= \mS_\CN \t \mS_{C'\CN}\t K(x_n)\ar[d]^v\\
\mS_{\CN}\t L_{C\CN}=\mS_\CN \t L_{C'\CN}\t K[x_n, x_n^{-1}]
\ar[r]^{\psi} &\mS_\CN \t L_{C'\CN}\t K(x_n)}$$ where the
horizontal maps are localization maps (they are injective) and the
vertical maps are epimorphisms (the maps $u$ and $v$ are
determined in the obvious way by the epimorphism $\pi_{C'\CN}:
\mS_{C'\CN}\ra L_{C'\CN}$).  Since $\overline{\gp}\cap
S=\emptyset$ and the set $S$ is central,
$$ \overline{\gp}= \v^{-1} v^{-1} (\mS_\CN \t \gq ) =u^{-1}
\psi^{-1} (\mS_\CN \t \gq ) = u^{-1} (\mS_\CN \t \psi_1^{-1} (\gq
))$$ where $\psi_1 : L_{C\CN}= L_{C'\CN}\t K[x_n,x_n^{-1}]\ra
L_{C'\CN}\t K(x_n)$ is the localization map, i.e. the restriction
of the map $\psi$ to the subalgebra $L_{C\CN}$ of $\mS_\CN \t
L_{C\CN}$. The ideal $\gq':= \psi_1^{-1}(\gq )$ is a prime ideal
of the algebra $L_{C\CN}$. Now, $\gp = \mS_\CN \t
\pi_{C\CN}^{-1}(\gq')$. The proof of the theorem is complete.
$\Box $

\begin{theorem}\label{d21Dec8}
The map $\Max (L_n) \ra \Max (\mS_n)$, $ \gm \mapsto \pi^{-1}(\gm
)$, is a bijection where $\pi : \mS_n\ra \mS_n/ \ga_n \simeq L_n$,
$a\mapsto a+\ga_n$, i.e. $\Max (\mS_n ) = \Spec (\mS_n, \emptyset
) = \{ \pi^{-1}(\gm )\, | \, \gm \in \Max (L_n)\}$.
\end{theorem}

{\it Proof}. The result follows from Theorem \ref{21Dec8} and the
fact that  the algebras $\mS_i$, $i\geq 1$, are not simple.
 $\Box $

$\noindent $

A prime ideal $\gp $ of an algebra $A$ is called a {\em completely
prime} ideal if the factor algebra $A/\gp $ is a domain.

\begin{corollary}\label{c9Jan9}
A prime ideal $\gp$ of the algebra $\mS_n$ is completely prime iff
$\ga_n\subseteq \gp$.
\end{corollary}

{\it Proof}. By Theorem \ref{21Dec8}.(2), $\gp = \mS_{\CN} \t
\pi^{-1}_{C\CN}(\gq )$. The factor algebra $\mS_n / \gp \simeq
\mS_\CN \t (\mS_{C\CN}/ \pi^{-1}_{C\CN}(\gq ))$ is a domain iff
$\CN = \emptyset$,  i.e. $\ga_n\subseteq \gp$.  $\Box $

\begin{corollary}\label{a9Jan9}
$\bigcap_{\gm \in \Max (\mS_n )}\gm = \ga_n$.
\end{corollary}

{\it Proof}. The statement follows at once from Theorem
\ref{d21Dec8} and the fact that $\bigcap_{\gm' \in \Max (L_n
)}\gm' =0$. $\Box $

$\noindent $

The next corollary gives all the primes of $\mS_n$ that are
contained in the prime ideal $\ga_n$. It turns out they are
precisely the primes of the type $\mS_\CN \t \pi^{-1}_{C\CN}(0)$
for all possible subsets $\CN$ of $\{ 1, \ldots , n\}$.

\begin{corollary}\label{g21Dec8}
\begin{enumerate}
\item The set $\{ \gp_\CN :=\sum_{i\in \CN}\gp_i\, | \, \CN
\subseteq \{ 1, \ldots , n\} \}$ is the set of all the prime
ideals of the algebra $\mS_n$ contained in the prime ideal
$\ga_n$, it contains precisely $2^n$ elements; $\gp_\CN
=\mS_{C\CN} \t \pi^{-1}_\CN (0)$; $\gp_\emptyset := 0$; $\gp_{\CN}
= \gp_\CM$ iff $\CN = \CM$. \item $\gp_\CN^2= \gp_\CN$ for all
$\CN$. \item $\gp_\CN \subseteq \gp_\CM$ iff $\CN \subseteq \CM$.
\item If $\gp_\CN \subseteq \gp_\CM$ then $\gp_\CN \gp_\CM=
\gp_\CN$. \item $\hht (\gp_\CN) = |\CN |$ for all $\CN$. In
particular, $\hht (\ga_n) = n$.
\end{enumerate}
\end{corollary}

{\it Proof}. 1. By Theorem \ref{21Dec8}, each prime ideal that is
contained in the ideal $\ga_n$ is of the type $\mS_\CN \t
\pi^{-1}_{C\CN}(0)=\mS_\CN \t \ga_{C\CN}= \gp_{C\CN}$, and vice
versa. By the very definition, $\gp_\CN = \gp_\CM$ iff $\CN = \CM$
($\gp_\CN = \mS_{C\CN}\t \pi^{-1}_\CN (0)$ and  $\gp_\CM =
\mS_{C\CM}\t \pi^{-1}_\CM (0)$).

2. $\gp_\CN = \sum_{i\in \CN}\gp_i= \sum_{i\in
\CN}\gp_i^2\subseteq (\sum_{i\in \CN}\gp_i)(\sum_{j\in \CN}\gp_j)=
\gp_\CN^2\subseteq \gp_\CN$, hence $\gp_\CN = \gp^2_\CN$.

3. Statement 3 is obvious.

4. $\gp_\CN \supseteq \gp_\CN \gp_\CM \supseteq \gp_\CN \gp_\CN =
\gp_\CN$, and so $\gp_\CN \gp_\CM = \gp_\CN$.

5. By statements 1 and  3, $\hht (\gp_\CN) \leq |\CN |$. Let $\CN
= \{ i_1, \ldots , i_s\}$. Then $0\subset \gp_{i_1}\subset \gp_{\{
i_1, i_2\}}\subset \cdots \subset \gp_\CN$ is the strictly
descending chain of primes of length $|\CN |$, and so $\hht
(\gp_\CN )= |\CN |$.  $\Box $

$\noindent $

An ideal $I$ of an algebra $A$ is called an {\em idempotent} ideal
if $I^2=I$. The next corollary gives all the idempotent prime
ideals of the algebra $\mS_n$. All the idempotent ideals of the
algebra $\mS_n$ are found (Theorem \ref{9Jan9}) and their
properties are studied in Section \ref{IDIDS}.

\begin{corollary}\label{b9Jan9}
Let $\gp$ be a prime ideal of the algebra $\mS_n$. The following
are equivalent.
\begin{enumerate}
\item $\gp$ is an idempotent ideal; \item  the ideal $\gp$ is
contained in all the maximal ideals of the algebra $\mS_n$; \item
the ideal $\gp$ is a sum of height one prime ideals of the algebra
$\mS_n$, i.e. $\gp = \gp_{\CN}$ (Corollary \ref{g21Dec8}).
\end{enumerate}

So, the set of all the idempotent prime ideals $\{ \gp_{\CN}  \, |
\, \CN \subseteq \{ 1, \ldots , n\} \}$ of the algebra $\mS_n$
contains precisely $2^n$ elements and its Krull dimension is $n$.
\end{corollary}

{\it Proof}. $(1\Rightarrow 2)$ Let $\gp$ be an idempotent prime
ideal of the algebra $\mS_n$. Then its image $\gp'$ under the
epimorphism $\mS_n\ra \mS_n/ \ga_n \simeq L_n$ is an idempotent
ideal, hence either $\gp'=0$ or $\gp'=L_n$. The second case is
impossible since the ideal $\ga_n$ is contained in all the maximal
ideals of the algebra $\mS_n$ (Corollary \ref{a9Jan9}), and so the
sum $\gp +\ga_n$ is contained in all the maximal ideals $\gm$ such
that $\gp \subseteq \gm$. This contradicts to the fact that $\gp
+\ga_n = \mS_n$.  Therefore, $\gp'\neq L_n$, and so $\gp\subseteq
\ga_n$. By Corollary \ref{a9Jan9}, the ideal $\gp$ is contained in
all the maximal ideals of the algebra $\mS_n$.

$(2\Rightarrow 3)$ See Corollary \ref{a9Jan9} and Corollary
\ref{g21Dec8}.

$(3\Rightarrow 1)$  See Corollary \ref{g21Dec8}.

The rest follows directly from Corollary \ref{g21Dec8}. $\Box $

$\noindent $

Let $\CN_1$ and $\CN_2$ be subsets of the set $\{ 1, \ldots , n\}$
such that $\CN_1\supseteq \CN_2$, and so $C\CN_1\subseteq C\CN_2$.
 Consider  the commutative diagram of  algebra homomorphisms:
$$\xymatrix{\mS_{C\CN_1}\ar[d]^{\pi_{C\CN_1}}\ar[r] & \mS_{C\CN_2}\ar[d]^{\pi_{C\CN_2}}\\
L_{C\CN_1}\ar[r] &L_{C\CN_2}}$$ where the horizontal maps are
natural monomorphisms. The next lemma gives necessary and
sufficient conditions for one prime ideal to contain another prime
ideal. It is instrumental in the proofs of Theorems \ref{f21Dec8}
and \ref{22Dec8}.

\begin{lemma}\label{e21Dec8}
Let $\gp_1' = \mS_{\CN_1}\t \pi^{-1}_{C\CN_1}(\gq_1)$ and  $\gp_2'
= \mS_{\CN_2}\t \pi^{-1}_{C\CN_2}(\gq_2)$ be prime ideals of the
algebra $\mS_n$. Then
\begin{enumerate}
\item $\gp_1'\subseteq \gp_2'$ iff $C\CN_1\subseteq C\CN_2$ and
$\gq_1\subseteq \gq_2$ (recall that $L_{C\CN_1}\subseteq
L_{C\CN_2}$). \item $\gp_1'= \gp_2'$ iff $C\CN_1=C\CN_2$ and
$\gq_1=\gq_2$. \item $\gp_1'\subset \gp_2'$ iff either
$C\CN_1\subset C\CN_2$ and $\gq_1\subseteq \gq_2$ or $C\CN_1=
C\CN_2$ and $\gq_1\subset \gq_2$.
\end{enumerate}
\end{lemma}

{\it Proof}. 1. $(\Rightarrow )$ Suppose that $\gp_1'\subseteq
\gp_2'$. Note that $\gp_{C\CN_1}=\mS_{\CN_1}\t
\pi^{-1}_{C\CN_1}(0)\subseteq \gp_1'$. Moreover,
$\gp_{C\CN_1}=\gp_1'\cap \ga_n $. Similarly,
$\gp_{C\CN_2}=\gp_2'\cap \ga_n $, and so $ \gp_{C\CN_1}\subseteq
\gp_{C\CN_2}$. By Corollary \ref{g21Dec8}.(3), $C\CN_1\subseteq
C\CN_2$. Then $L_{C\CN_1}\subseteq L_{C\CN_2}$ and
$$\ga_{C\CN_2}\subseteq \mS_{C\CN_2\backslash
C\CN_1}\t\pi^{-1}_{C\CN_1}(\gq_1 ) \subseteq
\pi^{-1}_{C\CN_2}(\gq_2)\subseteq \mS_{C\CN_2}.$$ Taking the
inclusions  modulo the ideal $\ga_{C\CN_2}$ of the algebra
$\mS_{C\CN_2}$, we have $L_{C\CN_2\backslash C\CN_1}\t
\gq_1\subseteq \gq_2\subseteq L_{C\CN_2}$ and so $\gq_1\subseteq
\gq_2$.

$(\Leftarrow )$ This implication is obvious.

2. $(\Rightarrow )$ If $\gp_1'=\gp_2'$ then
$$ \gp_{C\CN_1}= \gp_1'\cap\ga_n = \gp_2'\cap \ga_n =
\gp_{C\CN_2},$$ and so $C\CN_1=C\CN_2$ (Corollary
\ref{g21Dec8}.(1)), and finally $\CN_1= \CN_2$. Then
$\gq_1=\gq_2$.

$(\Leftarrow )$ This implication is obvious.

3. Statement 3 is an easy corollary of statements 1 and 2. $\Box $

\begin{theorem}\label{f21Dec8}
The classical Krull dimension  of the algebra $\mS_n$ is $2n$.
\end{theorem}

{\it Proof}. Since $\Kdim (L_{\CN})=|\CN |$ for all $\CN\subseteq
\{ 1, \ldots , n\}$, $\clKdim (\mS_n) \leq 2n$, by Theorem
\ref{21Dec8} and  Lemma \ref{e21Dec8}. In fact, $\clKdim (\mS_n) =
2n$ as the following strictly ascending chain of prime ideals of
length $2n$ shows:
\begin{eqnarray*}
 0&\subset & \gp_1 \subset \gp_1+\gp_2\subset \cdots \subset \ga_n:= \gp_1+\cdots + \gp_n \\
 &\subset & \ga_n+(x_1-1)\subset \cdots \subset \ga_n +
 (x_1-1, \ldots, x_n-1). \;\;\; \Box
\end{eqnarray*}

 For an ideal $\ga$ of a ring, let $\Min (\ga )$ denote the set of
 the minimal primes over $\ga$.

\begin{theorem}\label{25Feb9}
For each ideal $\ga$ of the algebra $\mS_n$, the set $\Min (\ga)$
of minimal primes over $\ga$ is a finite, non-empty set.
\end{theorem}

{\it Proof}. By Theorem \ref{f21Dec8}, $\Min (\ga ) \neq
\emptyset$. By Theorem \ref{21Dec8}.(2), it suffices to show that,
for each subset $\CN$ of the index set $\{ 1, \ldots , n\}$, there
are no more than finitely many (may be none) minimal primes over
$\ga$ of the type $\mS_\CN \t \pi^{-1}_{C\CN }(\gq )$ for some
prime ideal $\gq $ of the algebra $L_{C\CN }$. Suppose that this
is not the case. Then for some subset $\CN$ there are infinitely
many distinct minimal primes over $\ga$, say $\mS_\CN \t
\pi^{-1}_{C\CN }(\gq_i )$, $i\in \N$. Let $\{ \gq_1', \ldots ,
\gq_s'\}$ be the minimal primes of the ideal $\bigcap_{i\in
\N}\gq_i$ in the (commutative Noetherian) algebra $L_{C\CN}$. Then
$$ \ga \subseteq \mS_\CN \t \pi^{-1}_{C\CN }(\bigcap_{i\in \N}\gq_i )\subseteq \mS_\CN \t \pi^{-1}_{C\CN }(\gq_j'),
 \;\; j=1, \ldots , s, $$
and the ideals $\mS_\CN \t \pi^{-1}_{C\CN }(\gq_j')$ of the
algebra $\mS_n$ are prime. Then each ideal $\mS_\CN \t
\pi^{-1}_{C\CN }(\gq_i)$, $i\in \N$, must be equal to one of the
ideals $\mS_\CN \t \pi^{-1}_{C\CN }(\gq_j')$, a contradiction.
$\Box $

 $\noindent $

The next theorem provides formulae for the height and co-height of
primes of the algebra $\mS_n$. They are used in the proof of the
fact that the algebra $\mS_n$ is catenary (Theorem \ref{23Dec8}).

\begin{theorem}\label{22Dec8}
Let $\gp$ be a prime ideal of the algebra $\mS_n$, i.e. $\gp =
\mS_\CN \t \pi^{-1}_{C\CN}(\gq )$ (Theorem \ref{21Dec8}). Then
\begin{enumerate}
\item $\hht (\gp ) = |C\CN |+\hht_{L_{C\CN}}(\gq )$. \item $\cht
(\gp ) =2 |\CN |+\cht_{L_{C\CN}}(\gq )$.\item $\hht (\gp ) +\cht
(\gp ) = \clKdim (\mS_n)$.
\end{enumerate}
\end{theorem}

{\it Proof}. 1. By Lemma \ref{e21Dec8}.(3), $\hht (\gp ) \leq
|C\CN |+ \hht_{L_{C\CN}}(\gq )$. Note that the prime ideal
$\gp_{C\CN} = \mS_\CN \t \pi^{-1}_{C\CN}(0)$ has height $|C\CN |$
(Corollary \ref{g21Dec8}.(5)) and $\gp_{C\CN}\subseteq \gp$. The
length of the maximal chain of prime ideals of the algebra $\mS_n$
lying between the ideals $\gp_{C\CN}$ and $\gp$ is at least
$\hht_{L_{C\CN}}(\gq )$. Hence $\hht (\gp ) \geq |C\CN
|+\hht_{L_{C\CN}}(\gq )$. This proves the result.

2. \begin{eqnarray*}
 \cht (\gp ) &=&\clKdim (\mS_n / \gp ) = \clKdim (\mS_\CN \t \mS_{C\CN }/(\mS_\CN \t \pi^{-1}_{C\CN}(\gq ))) \\
 &=& \clKdim (\mS_\CN \t (L_{C\CN }/\gq ))\geq \clKdim (\mS_\CN )
 + \Kdim (L_{C\CN }/\gq ) \;\; ({\rm by\; Proposition}\;
 \ref{c21Dec8})\\
 &=& 2|\CN |+\cht_{L_{C\CN}}(\gq ) \;\; ({\rm by\; Theorem }\;
 \ref{f21Dec8}).
\end{eqnarray*}
Therefore, $\cht(\gp ) \geq  2|\CN |+\cht_{L_{C\CN}}(\gq )$. The
reverse inequality follows from the chain of inequalities below
which turn out to be equalities since the both ends are identical:
\begin{eqnarray*}
 \clKdim (\mS_n)&=&2n=2(|C\CN |+|\CN |) = |C\CN |+\Kdim (L_{C\CN})+2|\CN |\\
 &=& |C\CN |+\hht_{L_{C\CN}}(\gq) +\cht_{L_{C\CN}}(\gq ) +2|\CN |\\
 &=& \hht (\gp) + \cht_{L_{C\CN}}(\gq ) +2|\CN |\;\;\; ({\rm by \; statement}\; 1)\\
 &\leq & \hht (\gp ) +\cht (\gp ) \leq \clKdim (\mS_n).
\end{eqnarray*}

3. By statements 2 and 3,
\begin{eqnarray*}
\hht (\gp ) +\cht (\gp ) &=& |C\CN |+ \hht_{L_{C\CN}}(\gq ) +\cht_{L_{C\CN}}(\gq ) +2|\CN |\\
 &=& |C\CN |+|C\CN |+2|\CN |= 2(|C\CN |+|\CN |)=2n = \clKdim
 (\mS_n). \;\; \Box
\end{eqnarray*}

Let $\gp = \mS_\CN \t \pi^{-1}_{C\CN } (\gq )$ be a prime ideal of
the algebra $\mS_n$ (see Theorem \ref{21Dec8}). Let $s= |C\CN |$,
$C\CN = \{ i_1, \ldots , i_s\}$, $t=\hht_{L_{C\CN}}(\gq )$, and
$0\subset \gq_1\subset \gq_2\subset \cdots \subset \gq_t=\gq$ be a
strictly ascending chain of primes of the algebra $L_{C\CN }$.
Then the strictly ascending chain of primes in $\mS_n$ has length
$\hht (\gp )$:
\begin{eqnarray*}
 0&\subset &\gp_{i_1}\subset \gp_{i_1}+\gp_{i_2}\subset \cdots \subset \gp_{i_1}+\cdots + \gp_{i_s}=\mS_\CN \t \pi^{-1}_{C\CN}(0)
 \subset
 \mS_\CN \t \pi^{-1}_{C\CN}(\gq_1)\subset  \\
 &  \cdots & \subset \mS_\CN \t \pi^{-1}_{C\CN}(\gq_t)=\gp .
\end{eqnarray*}
Let $k= |\CN |$, $\CN = \{ j_1, \ldots , j_k\}$,  $l=
\cht_{L_{C\CN}}(\gq )$,  and $\gq \subset \gq_1'\subset
\gq_2'\subset \cdots \subset \gq_l'$ be a strictly ascending chain
of primes of the algebra $L_{C\CN }$. Let $I_i:= \mS_\CN \t
\pi^{-1}_{L_{C\CN}}(\gq_i')$. Then $I_0:= \gp \subset I_1\subset
\cdots \subset I_l$ is the strictly ascending chain of primes in
$\mS_n$. Consider the following strictly ascending chain of primes
in $\mS_n$ (see the proof of Theorem \ref{f21Dec8}):
\begin{eqnarray*}
 0&\subset & \gp_{j_1} \subset \gp_{j_1}+\gp_{j_2}\subset \cdots \subset \gp_\CN := \gp_{j_1}+\cdots + \gp_{j_k} \\
 &\subset & \gp_\CN+(x_{j_1}-1)\subset \cdots \subset \gp_\CN +
 (x_{j_1}-1, \ldots, x_{j_k}-1).
\end{eqnarray*}
Let us denote these ideals as $J_0:=0\subset J_1\subset \cdots
\subset J_{2k}$. Then the strictly ascending chain of primes in
$\mS_n$ has length $\cht (\gp )$:
$$ \gp = I_0+J_0\subset I_0+J_1\subset \cdots \subset
I_0+J_{2k}\subset I_1+J_{2k}\subset \cdots \subset I_l+J_{2k}.$$
In more detail, for each $s=1, \ldots , k$,
$$ \mS_n /(I_0+J_s)\simeq \mS_{\CN \backslash \{ j_1, \ldots ,
j_s\} } \t L_{\{ j_1, \ldots , j_s\} } \t (L_{C\CN } / \gq )$$ is
 a prime algebra, by Proposition \ref{c21Dec8}.(1) (since $ L_{\{ j_1, \ldots , j_s\} }\t (L_{C\CN } / \gq
 )$ is a prime algebra).

 For each $s= k+1, \ldots , 2k$, (i.e. $s= k+t$ where $t=1, \ldots
 , k$): $$\mS_n/(I_0+J_{k+t})\simeq L_{\{ j_{t+1}, \ldots , j_k\} }\t (L_{C\CN } / \gq
 )$$ is a prime algebra. For each $i=1, \ldots , l$,
 $\mS_n/(I_i+J_{2k})\simeq L_{C\CN }/ \gq_i'$ is a prime algebra.

$\noindent $

Let $A$ be an algebra and $\gp_1$, $\gp_2$ be its prime ideals
such that $\gp_1\subseteq \gp_2$. The {\em relative height} $\hht
(\gp_2, \gp_1)$ is the maximum of lengths of strictly ascending
chains of prime ideals of the type $\gp_1\subset \gq_2 \subset
\gq_3\subset \cdots \subset\gq_s \subset \gp_2$. If
$\gp_1\subseteq \gp_2\subseteq \gp_3$ then
$$\hht (\gp_3,
\gp_1)\geq \hht (\gp_3, \gp_2)+\hht (\gp_2, \gp_1).$$ Let
$\gp_1\subset \gp_2\subset\cdots \subset  \gp_s$ be a chain of
prime ideals of the algebra $A$. Consider a second chain of prime
ideals $\gp_1\subset \gp_{i_1}\subset\gp_{i_2}\subset \cdots
\subset\gp_{i_t}\subset  \gp_s$ where $1<i_1<i_2<\cdots < i_t<s$
which is obtained from the first chain by possibly deleting some
of its intermediate terms. We say that the first chain is a {\em
refinement} of the second chain.

$\noindent $

{\it Definition}. An algebra $A$ of finite classical Krull
dimension is called a {\em catenary} algebra if each chain of
prime ideals $\gp_1\subset \gp_2\subset\cdots \subset  \gp_s$
admits a refinement of length $\hht (\gp_s, \gp_1)$.

$\noindent $

The algebras $P_n$ and $L_n$ are catenary.

\begin{proposition}\label{a23Dec8}
Let $\gp_i=\mS_{\CN_i}\t \pi^{-1}_{C\CN_i}(\gq_i)$, $i=1,2,3$, be
prime ideals of the algebra $\mS_n$ such that $\gp_1\subseteq
\gp_2 \subseteq \gp_3$. Then
\begin{enumerate}
\item $\hht (\gp_2, \gp_1) = \hht (\gp_2) -\hht (\gp_1) = \cht
(\gp_1) - \cht (\gp_2)=
|\CN_1|-|\CN_2|+\hht_{L_{C\CN_2}}(\gq_2)-\hht_{L_{C\CN_1}}(\gq_1)$.
\item (Additivity of the relative height function): $\hht (\gp_3,
\gp_1)=\hht (\gp_3, \gp_2)+\hht (\gp_2, \gp_1)$.
\end{enumerate}
\end{proposition}

{\it Proof}. 1. By Theorem \ref{22Dec8}.(3), $\hht (\gp_1) +\cht
(\gp_1)=\hht (\gp_2) +\cht (\gp_2)$, and so the second equality
follows. The third equality follows from Theorem \ref{22Dec8}.(1),
\begin{eqnarray*}
 \hht (\gp_2) -\hht (\gp_1) & =&
 |C\CN_2|+\hht_{L_{C\CN_2}}(\gq_2)-|C\CN_1|-\hht_{L_{C\CN_1}}(\gq_1)\\
 &=&|\CN_1|-|\CN_2|+\hht_{L_{C\CN_2}}(\gq_2)-\hht_{L_{C\CN_1}}(\gq_1).
\end{eqnarray*}
Note that $\hht (\gp_2)= \hht (\gp_2, 0)\geq \hht (\gp_2, \gp_1)
+\hht (\gp_1, 0)= \hht (\gp_2, \gp_1)+\hht (\gp_1)$, and so $\hht
(\gp_2, \gp_1) \leq \hht (\gp_2) - \hht (\gp_1)$. It remains to
prove the inverse inequality, i.e.
$$\hht (\gp_2, \gp_1) \geq
|C\CN_2\backslash
C\CN_1|+\hht_{L_{C\CN_2}}(\gq_2)-\hht_{L_{C\CN_1}}(\gq_1).$$
First, we consider two special cases, and then the general case
can be deduced from these two special ones.

Case 1: $C\CN_1 = C\CN_2$.  Then $\gp_1= \mS_{\CN_1}\t
\pi^{-1}_{C\CN_1}(\gq_1)$, $\gp_2= \mS_{\CN_2}\t
\pi^{-1}_{C\CN_2}(\gq_2)$, and $\gq_1\subseteq \gq_2$. We may
assume that $\gq_1\neq \gq_2$. If $\gq_1\subset \gq_2'\subset
\cdots \subset \gq_s'\subset \gq_2$ is a strictly ascending chain
of primes of the algebra $L_{C\CN_1}$ then the induced chain
$\gp_1\subset \gp_2'\subset \cdots \subset \gp_s'\subset \gp_2$
(where $\gp_i':= \mS_{\CN_1}\t \pi^{-1}_{C\CN_1}(\gq_i')$) is a
strictly ascending chain of primes of the algebra $\mS_n$, and so
$$ \hht (\gp_2, \gp_1) \geq \hht_{L_{C\CN_1}}(\gq_2, \gq_1)=
\hht_{L_{C\CN_1}}(\gq_2) - \hht_{L_{C\CN_1}}(\gq_1), $$ as
required.

Case 2: $C\CN_1\subseteq C\CN_2$ and $L_{C\CN_2} \gq_1=\gq_2$.
Then $\CN_1\supseteq \CN_2$. We may assume that  $\CN_1\neq
\CN_2$. Fix a chain of sets $\CN_1\supset \CM_2\supset \cdots
\supset \CM_t\supset \CN_2$ such that the number of elements of
each successor is one less than of its predecessor. By Lemma
\ref{e21Dec8}.(3), we have the strictly ascending chain of primes
in $\mS_n$:
$$ \gp_1\subset \gp_2'\subset \cdots \subset \gp_t'\subset \gp_2$$
where $\gp_i':= \mS_{\CM_i} \t \pi^{-1}_{C\CM_i} (L_{C\CM_i}
\gq_1)$, and so
$$ \hht (\gp_2, \gp_1) \geq |C\CN_2\backslash C\CN_1|= |C\CN_2\backslash
C\CN_1|+\hht{L_{C\CN_2}}(L_{C\CN_2}\gq_1) -
\hht{L_{C\CN_1}}(\gq_1), $$ as required since
$\hht_{L_{C\CN_2}}(L_{C\CN_2}\gq_1) = \hht_{L_{C\CN_1}}(\gq_1)$.

In general case, we have the inclusions of the prime ideals (Lemma
\ref{e21Dec8}):
$$ \gp_1= \mS_{\CN_1} \t \pi^{-1}_{C\CN_1}(\gq_1)\subseteq \gp :=
\mS_{\CN_2}\t \pi^{-1}_{C\CN_2}(L_{C\CN_2}\gq_1) \subseteq \gp_2=
\mS_{\CN_2} \t \pi^{-1}_{C\CN_2}(\gq_2), $$ where the first
inclusion is Case 2, and the second inclusion is Case 1.

 Now,
\begin{eqnarray*}
 \hht (\gp_2, \gp_1)&\geq & \hht (\gp_2, \gp )+\hht (\gp , \gp_1)\geq  \hht_{L_{C\CN_2}}(\gq_2) -
\hht_{L_{C\CN_2}}(L_{C\CN_2}\gq_1)+|C\CN_2\backslash C\CN_1|\\
 &=&\hht_{L_{C\CN_2}}(\gq_2) -
\hht_{L_{C\CN_1}}(\gq_1)+|C\CN_2\backslash C\CN_1|,
\end{eqnarray*}
as required.

2. Statement 2 follows at once from statement 1.  $\Box $

\begin{theorem}\label{23Dec8}
The algebra $\mS_n$ is a catenary algebra.
\end{theorem}

{\it Proof}. The fact that the algebra $\mS_n$ is a catenary
algebra follows at once from the additivity of the relative height
function (Proposition \ref{a23Dec8}.(2)). In more detail, given a
chain of prime ideals $\gp_1\subset \gp_2\subset \cdots \subset
\gp_s$. Then $\hht (\gp_s, \gp_1) = \sum_{i=2}^s\hht (\gp_i,
\gp_{i-1})$. So, the algebra $\mS_n$ is catenary.  $\Box $

$\noindent $


\section{Left or right Noetherian factor algebras of $\mS_n$}

The aim of this short  section is to give answers (Theorem
\ref{24Dec8}) to the following two questions:

$\noindent $

{\it Question 1. For which ideals $I$ of $\mS_n$, the factor
algebra $\mS_n/ I$ is a left or right Noetherian algebra?}

$\noindent $

{\it Question 2. For which ideals $I$ of $\mS_n$, the factor
algebra $\mS_n/ I$ is a commutative algebra?}

$\noindent $

{\it Example}. $I=\ga_n$ since $\mS_n / \ga_n\simeq L_n$ is a
commutative Noetherian algebra.
\begin{theorem}\label{24Dec8}
Let $I$ be an ideal of the algebra $\mS_n$. The following are
equivalent.
\begin{enumerate}
\item The factor algebra $\mS_n/I$ is a left Noetherian algebra.
\item The factor algebra $\mS_n/I$ is a right  Noetherian algebra.
\item The factor algebra $\mS_n/I$ is a commutative algebra. \item
$\ga_n\subseteq I$.
\end{enumerate}
\end{theorem}

{\it Proof}. It suffices to show that $1\Leftrightarrow 3
\Leftrightarrow 4$ since then applying the involution $\eta$ we
have the equivalences $ 2 \Leftrightarrow 4$ since $\eta (\ga_n) =
\ga_n$.

$(4\Rightarrow 3)$ The algebra $\mS_n/I$ is commutative as a
factor algebra of the commutative algebra $L_n\simeq  \mS_n/
\ga_n$.

$(3\Rightarrow 1)$ Trivial.

$(1\Rightarrow 4)$ Suppose that the factor algebra $\mS_n/I$ is
left Noetherian. We have to show that $\ga_n\subseteq I$.  Since
$\ga_n = \gp_1+\cdots +\gp_n$, we have to show that all
$\gp_i\subseteq I$. By symmetry of the indices $1, \ldots , n$, it
suffices to show that $\gp_1\subseteq I$. By Lemma
\ref{c23Dec8}.(1), $I\cap \gp_1= F\t I_{n-1}$ for some ideal
$I_{n-1}$ of the algebra $\mS_{n-1}$. Then
$$ \mS_n / I\supset \gp_1/ (I\cap \gp_1) = F\t \mS_{n-1}/(F\t
I_{n-1}) \simeq F\t (\mS_{n-1}/ I_{n-1}),$$ and so $F\t
(\mS_{n-1}/ I_{n-1})$ is a Noetherian $\mS_n$-module. Since the
$\mS_1$-module $F$  is an infinite  direct sum of nonzero
submodules, we deduce that $I_{n-1}= \mS_{n-1}$. Then $I\cap
\gp_1= F\t \mS_{n-1}= \gp_1$, and so $\gp_1\subseteq I$, as
required.  $\Box $

$\noindent $


\section{The weak and  global dimensions of the algebra $\mS_n$}

In this section, it is shown that the weak dimension and the left
and right global  dimensions of the algebra $\mS_n$ are equal to
$n$ (Theorem \ref{A11Dec8} and Corollary \ref{26Dec8}); the
projective dimensions of all simple $\mS_n$-modules $M$ are found
explicitly (Theorem \ref{28Dec8}), and it is shown that $\pd (M)
+\GK (M) = \lgldim (\mS_n)$ (Corollary \ref{x28Dec8}).

It follows from the decomposition of the  vector space $\mS_n =
P_n\oplus(\sum_{i=1}^n \mS_n y_i)$ that 
\begin{equation}\label{SnPn}
{}_{\mS_n}P_n\simeq \mS_n / (\sum_{i=1}^n \mS_n y_i).
\end{equation}
\begin{proposition}\label{a11Dec8}
\begin{enumerate}
\item The left $\mS_n$-module $P_n$ is projective.\item
$F_n:=F^{\t n }$ is a left and right projective $\mS_n$-module.
\item The projective dimension of the left and right
$\mS_n$-module $\mS_n/F_n$ is 1. \item For each element $\alpha
\in \N^n$, the $\mS_n$-module $\mS_n/ \mS_n y^\alpha$ is
projective.
\end{enumerate}
\end{proposition}

{\it Proof}. 1. The fact that the left $\mS_n$-module $P_n$ is
projective follows at once from the decomposition 
\begin{equation}\label{S1y}
\mS_1 = \mS_1y \oplus \mS_1E_{00},
\end{equation}
the fact that $\mS_n = \mS_1^{\t n }$ and that
${}_{\mS_1}\mS_1E_{00}\simeq K[x]$. In more detail, by
(\ref{S1y}),
\begin{eqnarray*}
 \mS_n&=& \mS_1(1)\t \cdots \t \mS_1(n)=\bigotimes_{i=1}^n (\mS_1(i)y_i\oplus \mS_1 (i) E_{00}(i))
 = \mS_n E_{00}(1)\t\cdots \t E_{00}(n)\bigoplus (\sum_{i=1}^n \mS_n y_i)\\
 &\simeq & P_n \bigoplus (\sum_{i=1}^n \mS_n y_i)
\end{eqnarray*}
since ${}_{\mS_n}P_n\simeq \mS_n E_{00}(1)\t\cdots \t E_{00}(n)$,
and so $P_n$ is a projective $\mS_n$-module.

 In order to prove (\ref{S1y}),
let us consider the short exact sequence of $\mS_1$-modules:
$$ 0\ra \mS_1\stackrel{r_y}{\ra}\mS_1\ra
\mS_1/\mS_1y\simeq K[x]\ra 0$$ where $r_y(a)=ay$ (note that
${}_{\mS_1}\mS_1y\simeq \mS_1$, by Lemma \ref{a7Dec8}.(2)). The
homomorphism $r_y$ admits a splitting given by the homomorphism
$r_x:\mS_1\ra \mS_1$, i.e. $r_xr_y=\id_{\mS_1}$ since $yx=1$. Now,
$\mS_1= \mS_1y \bigoplus \im ({\rm id}_{\mS_1} - r_yr_x)$. Since
${\rm id}_{\mS_1} - r_yr_x=r_{1-xy}= r_{E_{00}}$, the
decomposition (\ref{S1y}) becomes obvious.

2. Note that $F=\bigoplus_{i\geq 0} \mS_1E_{ii}\simeq
\bigoplus_{i\geq 0}P_1$ (the $\mS_1$-module $F$ is a direct sum of
infinite number of copies of the $\mS_1$-module $P_1$). Now,
\begin{equation}\label{SnFn}
{}_{\mS_n}F_n\simeq F^{\t n}\simeq \bigoplus_{\alpha\in \N^n}P_n.
\end{equation}
By statement 1, the left $\mS_n$-module $F_n$ is projective. Using
the involution $\eta$ and the fact that $\eta (F_n) = F_n$, we see
that the right $\mS_n$-module $F_n$ is projective.

3. The short exact sequence of left and right $\mS_n$-modules
$0\ra F_n \ra \mS_n\ra \mS_n / F_n\ra 0$ does not split since
$F_n$ is an essential left and right submodule of $\mS_n$
(Proposition \ref{b21Dec8}.(3)). By statement 1, the projective
dimension of the left and right $\mS_n$-module $\mS_n / F_n$ is 1.

4. Let $\Z^n =\bigoplus_{i=1}^n \Z e_i$ where $e_1, \ldots , e_n$
is the canonical free $\Z$-basis for $\Z^n$. Let $m=|\alpha |$.
Fix a chain of elements of $\Z^n$, $\beta_0=0, \beta_1, \ldots ,
\beta_m=\alpha$ such that, for all $i$, $\beta_{i+1}=\beta_i+e_j$
for some index $j=j(i)$. Then all the factors of the chain of left
ideals
$$\mS_ny^\alpha =\mS_ny^{\beta_m} \subset \mS_ny^{\beta_{m-1}} \subset\cdots \subset \mS_ny^{\beta_1}
\subset \mS_n$$ are projective $\mS_n$-modules since
$\mS_ny^{\beta_i}/\mS_ny^{\beta_{i+1}}\simeq \mS_n/\mS_ny_j\simeq
K[x_j]\t\mS_{n-1}$ is the projective $\mS_n$-module (statement 1).
Therefore, the $\mS_n$-module $\mS_n/ \mS_ny^\alpha $ is
projective. $\Box $

$\noindent $

We are interested in the homological properties of the algebra
$L_n\simeq \mS_n / \ga_n$ as a left and right $\mS_n$-module. For
a ring $R$ and modules ${}_RM$ and $N_R$, we denote by
$\pd({}_RM)$ and $\pd(N_R)$ their projective dimensions.

\begin{proposition}\label{b11Dec8}
$\pd ({}_{\mS_1}L_1) = \pd ({L_1}_{\mS_1})=1$.
\end{proposition}

{\it Proof}. Due to the involution $\eta$ of the algebra $\mS_1$
and the fact that $\eta (F) = F$ it suffices to prove only the
first statement, i.e. $\pd ({}_{\mS_1}L_1) =1$.

 Recall that the left $\,S_1$-module $F$ is
projective (Proposition \ref{a11Dec8}). So, the short exact
sequence of $\mS_1$-modules 
\begin{equation}\label{F1L1}
0\ra F\ra \mS_1\ra \mS_1/ F\simeq L_1\ra 0
\end{equation}
is a projective resolution for the $\mS_1$-module $L_1$. It
suffices to prove that the short exact sequence does not split
(this fact then implies that $\pd ({}_{\mS_1}L_1) =1$). Suppose
that this is not the case. Then there exists a homomorphism
$\mS_1\ra F$, $1\mapsto f$, such that $af=a$ for all elements
$a\in F$, or, equivalently, $F(f-1)=0$. By Proposition
\ref{10Dec8}.(2), $f=1\in F$, a contradiction. Therefore, $\pd
({}_{\mS_1}L_1) =1$.  $\Box $

$\noindent $

A module is called a {\em cyclic} module if it is generated by a
single element.

\begin{lemma}\label{c11Dec8}
The projective dimension of all the  nonzero cyclic left or right
$\mS_1$-modules annihilated by the ideal $F$ is 1.
\end{lemma}

{\it Proof}. In view of the involution $\eta$ of the algebra
$\mS_1$ it suffices to prove the result, say, for left modules
since $\eta (F) = F$.

Let $M$ be a nonzero cyclic $\mS_1$-module which is annihilated by
the ideal $F$. This means that $M$ is a cyclic $L_1$-module since
$\mS_1/F\simeq L_1=K[x,x^{-1}]$, and so either $M\simeq L_1$ or,
otherwise, $M\simeq M_a:= L_1/L_1a$ for some non-scalar polynomial
$a=a(x)\in K[x]$ such that $a(0)\neq 0$ (Proposition
\ref{10Dec8}.(4)). By Proposition \ref{b11Dec8},
$\pd_{\mS_1}(L_1)=1$. So, let $M=M_a$. We claim that
\begin{equation}\label{lan}
\lann_{\mS_1}(a):=\{ b\in \mS_1\, | \, ba=0\}=0.
\end{equation}
Let $I:=\lann_{\mS_1}(a)$. Since the element $a$ is a regular
element of the algebra $L_1$, it follows from the short exact
sequence of $\mS_1$-modules: $0\ra F\ra \mS_1\ra L_1\ra 0$ that
$I\subseteq F$. Since $E_{ij} x = E_{i,j-1}$ (see (\ref{xyEij1})),
$a(0)\neq 0$, and $F=\bigoplus_{i,j\geq 0}KE_{ij}$, we see that
$I=0$, as required. By (\ref{lan}), ${}_{\mS_1}\mS_1a\simeq \mS_1$
and ${}_{\mS_1}F\simeq Fa$. Since ${}_{\mS_1}F\simeq Fa\subseteq
F\cap \mS_1a\subseteq F\simeq P_1^{(\N )}$, we see that
${}_{\mS_1}F\cap \mS_1a\simeq P_1^{(\N )}\simeq F$. Using this
fact we can produce a splitting to the homomorphism $\alpha$ in
the short exact sequence of $\mS_1$-modules:
$$ 0\ra F\cap \mS_1a\stackrel{\alpha}{\ra}F\oplus
\mS_1a\stackrel{\beta}{\ra} F+\mS_1 a\ra 0, \;\; \alpha (u) =
(u,-u), \;\; \beta (u,v) = u+v.$$ Therefore, the $\mS_1$-modules
$F+\mS_1a$ is projective, and so the short exact sequence of
$\mS_1$-modules: 
\begin{equation}\label{Ma}
0\ra F+ \mS_1 a\ra \mS_1\ra \mS_1/ (F+\mS_1a)\simeq M_a\ra 0
\end{equation}
is a projective resolution for the $\mS_1$-module $M_a$. This
sequence is not split since the algebra $\mS_1$ does not contain
finite dimensional left ideals. Therefore, $\pd_{\mS_1}(M_a)=1$.
$\Box $

$\noindent $

Let $\lgldim $ and $\rgldim$ stand for the left and right global
dimension of algebra respectively.

\begin{theorem}\label{11Dec8}
$\lgldim (\mS_1)= \rgldim (\mS_1)=1$.
\end{theorem}

{\it Proof}. The algebra $\mS_1$ is self-dual, and so $\lgldim
(\mS_1)= \rgldim (\mS_1)$. Let us prove that $\lgldim (\mS_1)=1$.
Note that
$$ \lgldim (\mS_1)= \max \{ \pd (\mS_1/I)\, | \, 0 \neq I\subseteq
{}_{\mS_1}\mS_1\}.$$ Using the $3\times 3$ Lemma, the short exact
sequence of $\mS_1$-modules: $0\ra F\ra \mS_1\ra L_1\ra 0$ and the
inclusion $I\subseteq \mS_1$ yield the short exact sequence of
$\mS_1$-modules:
$$0\ra F/(F\cap I)\ra \mS_1/I\ra \mS_1/(F+I)\simeq L_1/I'\ra 0$$
for some ideal $I'$ of the algebra $L_1$. The $\mS_1$-module
$F/(F\cap I)$ is either a zero one or, otherwise, is isomorphic to
a direct sum of several copies (may be infinitely many) of the
projective $\mS_1$-module $K[x]$ (Proposition \ref{a11Dec8}),
hence $\pd (F/(F\cap I))=0$, in this case.  Then $\pd (\mS_1/I)
\leq  \pd (L_1/I')\leq 1$ (Lemma \ref{c11Dec8}). Therefore,
$\lgldim (\mS_1)=1$.  $\Box $

$\noindent $

{\bf The  weak and global dimensions of the algebra $\mS_n$}. Let
$S$ be a non-empty multiplicatively closed subset of a ring $R$,
and let $\ass (S):= \{ r\in R\, | \, sr =0$ for some $s\in S\}$.
Then a {\em left quotient ring} of $R$ with respect to $S$ is a
ring $Q$ together with a homomorphism $\v : R\ra Q$ such that

(i) for all $s\in S$, $\v (s)$ is a unit in $Q$;

(ii) for all $q\in Q$, $q= \v (s)^{-1} \v (r)$ for some $r\in R$
and $s\in S$, and

(iii) $\ker (\v ) = \ass (S)$.

If there exists a left quotient ring $Q$ of $R$ with respect to
$S$ then it is unique up to isomorphism, and it is denoted
$S^{-1}R$. It is also said that the ring $Q$ is the {\em left
localization} of the ring $R$ at $S$.

{\it Example 1}. Let $S:= S_y:= \{ y^i, i\geq 0\}$ and $R=\mS_1$.
By (\ref{xyEij}), $\ass (S) =F$, and,  by (\ref{mS1d1}),
$\mS_1/\ass (S)\simeq K[y,y^{-1}]$. This means that the ring
$L_1=\mS_1/F$ is the left quotient ring of $\mS_1$ at $S$.

{\it Example 2}. Let $S_{y_1, \ldots , y_n}:= \{ y^\alpha , \alpha
\in \N^n\}$ and $R=\mS_n$. Then $\ass (S_{y_1, \ldots , y_n}) =
\ga_n$ and $\mS_n / \ga_n \simeq L_n$. Therefore,
\begin{equation}\label{Sy1yn}
S_{y_1, \ldots , y_n}^{-1}\mS_n\simeq L_n,
\end{equation}
i.e. $L_n$ is the left quotient ring of $\mS_n$ at $S_{y_1, \ldots
, y_n}$. Note that the right localization $\mS_n S_{y_1, \ldots ,
y_n}^{-1}$ of $\mS_n$ at $S_{y_1, \ldots , y_n}$ does not exist.
Otherwise, we would have $S_{y_1, \ldots , y_n}^{-1}\mS_n\simeq
\mS_n S_{y_1, \ldots , y_n}^{-1}$ but all the elements $y^\alpha$
are left regular, and we would have a monomorphism $\mS_n\ra
S_{y_1, \ldots , y_n}^{-1}\mS_n\simeq L_n$, which would be
impossible since the algebra $L_n$ is commutative but the algebra
$\mS_n$ is not.

$\noindent $

Let us recall a result which will be used repeatedly in the proof
of Theorem \ref{A11Dec8}.

\begin{proposition}\label{Aus-sub}
{\rm \cite{Aus-NagoyaMJ-55}}. Let $M$ be a module over an algebra
$A$, $I$ a non-empty well-ordered set, $\{ M_i\}_{i\in I}$ be  a
family of submodules of $M$ such that if $i,j\in I$ and $i\leq j$
then $M_i\subseteq M_j$. If $M=\bigcup_{i\in I}M_i$ and
$\pd_A(M_i/M_{<i})\leq n$ for all $i\in I$ where $M_{<i}:=
\bigcup_{j<i}M_j$ then $\pd_A(M)\leq n$.
\end{proposition}
Let $V\subseteq U\subseteq W$ be modules. Then the factor module
$U/V$ is called a {\em sub-factor} of the module $W$. Let $\wdim$
  denote  the  {\em weak global
dimension}.

\begin{theorem}\label{A11Dec8}
$\lgldim (\mS_n)= \rgldim (\mS_n)=n$.
\end{theorem}

{\it Proof}. In view of the involution $\eta$, $\lgldim (\mS_n) =
\rgldim (\mS_n)$. By (\ref{Sy1yn}),
\begin{equation}\label{vss3}
 n=\wdim (L_n) = \wdim (S_{y_1, \ldots , y_n}^{-1}\mS_n)\leq
\wdim (\mS_n) \leq \lgldim (\mS_n).
\end{equation}
 It remains to prove that $\lgldim (\mS_n) \leq
n$. Let $\bK$ be the algebraic  closure of the field $K$. Since
$\lgldim (\mS_n) \leq \lgldim (\bK \t \mS_n)$, we may assume that
$K=\bK$. To prove the inequality $\lgldim (\mS_n) \leq n$ we use
induction on $n$. The case $n=1$ has been considered above
(Theorem \ref{11Dec8}). Let $n>1$, and we assume that the result
holds for all $n'<n$. The algebra $K(y)\t\mS_{n-1}$ is the
localization of the algebra $K[y,y^{-1}]\t \mS_{n-1}$ at the
multiplicative set $K[y,y^{-1}]\backslash \{ 0\}$, hence the
algebra $K(y)\t\mS_{n-1}$ is the left localization $T^{-1}\mS_n$
of the algebra $\mS_n$ at the multiplicative set $T:=
K[y]\backslash \{ 0 \}$. By induction, $\lgldim
(T^{-1}\mS_n)=n-1$. Let $M$ be an $\mS_n$-module. We have to prove
that $\pd (M)\leq n$. To the localization map $\v : M\ra T^{-1}M$,
$ m\mapsto m/1$, we attach two short exact sequences of
$\mS_n$-modules: 
\begin{equation}\label{vss1}
0\ra M':= \ker (\v ) \ra M \ra \bM := \im (\v )\ra 0,
\end{equation}
\begin{equation}\label{vss2}
0\ra \bM \ra T^{-1}M\ra M'':=\coker (\v ) \ra 0.
\end{equation}
The $\mS_n$-module $M'=\{ m\in M\, | \,  tm=0$ for some element
$t\in T\}$ is the $T$-{\em torsion} submodule of $M$. To show that
$\pd (M)\leq n$ it suffices to show that the projective dimensions
of the $\mS_n$-modules $M'$, $T^{-1}M$, and $M''$ are less or
equal to $n$. Indeed, then by (\ref{vss2}), $\pd (\bM ) \leq \max
\{ \pd (T^{-1}M), \pd (M'')\} \leq n$; and,  by (\ref{vss1}), $\pd
(M) \leq \max \{ \pd (M'), \pd (\bM ) \} \leq n$.

An $\mS_n$-module $N$ is called a $T$-{\em torsion} $\mS_n$-module
if $T^{-1}N=0$. The modules $M'$ and $M''$ are $T$-torsion. The
fact that their projective dimensions are less or equal to $n$
follows from the claim below.

{\it Claim. $\pd (N)\leq n $ for all $T$-torsion $\mS_n$-modules
$N$}.

{\it Proof of the Claim}. By Theorem \ref{a17Dec8}.(1),
$$ \hmS_1= \{ K[x]\} \cup \widehat{K[y,y^{-1}]} = \{ K[x], K_\l := K[y,y^{-1}]/ (y-\l)\, | \,
\l \in K^*:=K\backslash \{ 0 \} \,
\}, $$ and so each simple $\mS_1$-module $U$ is $T$-torsion with
$\End_{\mS_1}(U)\simeq K$ (Lemma \ref{a17Dec8}.(2)) since $K=
\bK$. Any nonzero $T$-torsion $\mS_1$-module $V$ contains a
nonzero submodule $V'$ which is an epimorphic image of the
$\mS_1$-module $M_\l := \mS_1/ \mS_1(y-\l )$ for some element $\l
\in K$.

If $\l =0$ then $M_0=\mS_1/ \mS_1 y \simeq K[x]$ is a simple
$\mS_1$-module, and so is the $\mS_1$-module $V'\simeq K[x]$.

If $\l \neq 0$ then it follows from the decomposition $\mS_1=
\bigoplus_{i,j\geq 0} Kx^iy^j$ that $M_\l = \bigoplus_{i\geq 0}
Kx^i \overline{1}$ where $\overline{1}:= 1+\mS_1 (y -\l )$. For
all $i\geq 0$,
$$ E_{i0} \overline{1}=x^i E_{00}\overline{1}=x^i(1-xy) \overline{1}=x^i (1-\l x)
\overline{1}\neq 0, $$ and so ${}_{\mS_1}K[x]\simeq \mS_1E_{00}
\overline{1}=\bigoplus_{i\geq 0} Kx^i(1-\l x)\overline{1}$. There
is the short exact sequence of $\mS_1$-modules: 
\begin{equation}\label{MlKl}
0\ra K[x]\simeq \mS_1E_{00}\overline{1}\ra M_\l \ra K_\l \ra 0,
\end{equation}
which is non-split since ${}_{K[x]}M_\l \simeq K[x]$. Therefore,
the socle of the $\mS_1$-module $M_\l$ is $\mS_1E_{00}\simeq K[x]$
which is a simple, essential submodule of $M_\l$.
 Then either
$V'\simeq M_\l$ or, otherwise, $V'\simeq K_\l$. It follows that
the $\mS_1$-module $V'$ contains a simple, $T$-torsion
$\mS_1$-submodule. Any sub-factor of the $T$-torsion
$\mS_n$-module $N$ is a $T$-torsion $\mS_n$-module, and so it
contains a nonzero submodule of the type $U\t N'$ for some simple,
$T$-torsion $\mS_1$-module $U$ and some $\mS_{n-1}$-module $N'$.
By Proposition \ref{Aus-sub},
$$ \pd (N)\leq  \max \{ \pd (U\t N')\, | \, U\in \hmS_1,
{}_{\mS_{n-1}}N'\} \leq \lgldim (\mS_1) + \lgldim (\mS_{n-1})=
1+(n-1) =
 n. $$
The proof of the Claim is complete.

It remains to show that $\pd(T^{-1}M)\leq n$. By Theorem
\ref{11Dec8}, $\pd_{\mS_1}(K(y))\leq 1$ where $K(y) =
T^{-1}(\mS_1/ F)$. In fact, 
\begin{equation}\label{pdKy}
\pd_{\mS_1}(K(y))= 1,
\end{equation}
as it follows from the short exact sequence of $\mS_1$-modules:
$$0\ra K[y,y^{-1}]\simeq \mS_1/F\ra K(y) \ra K(y) / K[y,y^{-1}]\ra
0$$ and from Proposition \ref{b11Dec8} (if $\pd_{\mS_1}(K(y))=0$
then $2\leq \pd_{\mS_1}(K(y) / K[y,y^{-1}])\leq 1$, a
contradiction). Then
$$ \pd_{\mS_n}(T^{-1}\mS_n) = \pd_{\mS_1\t \mS_{n-1}}(K(y) \t
\mS_{n-1})\leq \pd_{\mS_1}(K(y))\leq 1,$$ and so the projective
dimension over the algebra $\mS_n$ of each projective
$T^{-1}\mS_n$-module does not exceed $1$. Using the Ext-long
sequence, we deduce that $\pd_{\mS_n} ( T^{-1}M) \leq \pd_{K(y) \t
\mS_{n-1}}(T^{-1}M) + 1 \leq (n-1)+1=n.$ The proof of the theorem
is complete. $\Box $

\begin{corollary}\label{aA11Dec8}
Suppose, in addition, that the field $K$ is algebraically closed.
Then each $T$-torsion $\mS_n$-module $N$ where $T=K[y]\backslash
\{ 0\}$ is a union $N=\cup_{i\in I} N_i$ of its submodules $N_i$
for some well-ordered set $I$ such that $i\leq j$ implies
$N_i\subseteq M_j$, and $N_i/N_{<i}\simeq U_i\t N_i'$ for all $i$
where $U_i$ is a simple,  $T$-torsion  $\mS_1$-module and $N_1'$
is an $\mS_{n-1}$-module.
\end{corollary}

{\it Proof}. The result follows at once from the fact established
in the proof of the Claim above saying  that any nonzero
sub-factor of a $T$-torsion $\mS_n$-module contains a submodule of
the type $U\t N'$ for some simple,  $T$-torsion  $\mS_1$-module
$U$ and an $\mS_{n-1}$-module $N'$. $\Box $


\begin{corollary}\label{26Dec8}
$\wdim (\mS_n)=n$.
\end{corollary}

{\it Proof}. The result follows from Theorem \ref{A11Dec8} and
(\ref{vss3}).  $\Box $

$\noindent $

{\bf The projective dimensions of simple $\mS_n$-modules}. Our aim
is to find the projective dimension of each simple $\mS_n$-module
(Theorem \ref{28Dec8}). For, we need the next lemma and corollary.

For each vector $\l = (\l_1, \ldots , \l_n) \in K^n$, consider the
$\mS_n$-module
$$ M_\l := \mS_n/\sum_{i=1}^n \mS_n (y_i-\l_i)\simeq
\bigotimes_{I=1}^n \mS_1(i)/\mS_1(i)(y_i-\l_i).$$

\begin{lemma}\label{a10Jan9}
For  each element  $\l = (\l_1, \ldots , \l_n) \in K^{*n}$,  the
projective dimension of the $\mS_n$-module $M_\l$ is $n$.
Moreover, $\Ext_{\mS_n}^n(M_\l , \mS_n)\simeq \mS_n /
\sum_{i=1}^n(y_i-\l_i) \mS_n\simeq K.$
\end{lemma}

{\it Proof}. The algebra $\mS_n$ is the bimodule
${}_{\mS_n}{\mS_n}_\CP$ where $\CP:= K[y_1-\l_1, \ldots ,
y_n-\l_n]$ is the polynomial subalgebra of $\mS_n$ in $n$
variables. The sequence $s=(y_1-\l_1, \ldots , y_n-\l_n)$ is a
{\em regular} sequence of the right $\CP$-module $\mS_n$. The
Koszul complex $K(s)$ yields a projective resolution for the right
$\CP$-module  $M_\l$. The differential of the Koszul complex
$K(s)$ is obviously an $\mS_n$-homomorphism. Therefore, $K(s)$ is
a projective resolution for the $\mS_n$-module $M_\l$. Using this
resolution it is easy to show that
$$\Ext_{\mS_n}^n(M_\l ,
\mS_n)\simeq \mS_n / \sum_{i=1}^n(y_i-\l_i) \mS_n\simeq
\bigotimes_{i=1}^n \mS_1(i)/(y_i-\l_i)\mS_1(i).$$ To finish the
proof it remains to show that $ \mS_1(i)/(y_i-\l_i)\mS_1(i)\simeq
K$ since then $\pd_{\mS_n}(M_\l ) =n$.

By symmetry of indices, it suffices to show that $ \mS_1/(y-\l
)\mS_1\simeq K$ for $\l\in K^*$. By (\ref{xyEij}), the linear map
$l_{y-\l }: F\ra F$, $f\mapsto (y-\l ) f$, is a bijection since
$\l \neq 0$. The Snake Lemma for the commutative diagram $$
\xymatrix{0\ar[r]&F\ar[d]_{l_{y-\l }}\ar[r] & \mS_1\ar[d]_{l_{y-\l }}\ar[r]& L_1\ar[d]_{l_{y-\l }} \ar[r] &0\\
0\ar[r] &F\ar[r] & \mS_1\ar[r] & L_1  \ar[r] & 0\\}
$$
gives a vector space isomorphism
$$ \mS_1/(y-\l ) \mS_1=\coker_{\mS_1}(l_{y-\l }) \simeq \coker_{L_1}(l_{y-\l
})= L_1/(y-\l )L_1 \simeq K,$$ as required. $\Box $

\begin{corollary}\label{b10Jan9}
For  each element  $\l = (\l_1, \ldots , \l_n) \in K^n$,  the
projective dimension of the $\mS_n$-module $M_\l$ is $n(\l
):=\#\{\l_i\, | \, \l_i\neq 0\}$. Moreover, $\Ext_{\mS_n}^{n(\l
)}(M_\l , \mS_n)\simeq P_{n-n(\l )}.$
\end{corollary}

{\it Proof}. Let $k=n(\l )$. If $k=n$ then the result is Lemma
\ref{a10Jan9}. So, we can assume that $k\neq n$. Then $M_\l \simeq
P_{n-k}\t M_\mu$ (up to order of the indices)  where the vector
$\mu = (\mu_1, \ldots , \mu_k) \in K^{*k}$ consists of all the
nonzero coordinates of the vector $\l \in K^n$. Now, applying
$P_{n-k}\t -$ to the Koszul complex $K(s')$ where $s'=(y_1-\mu_1,
\ldots , y_k-\mu_k)$ for the $\mS_k$-module $M_\mu$ (from the
proof of Lemma \ref{a10Jan9}) we obtain  a projective resolution
for the $\mS_n$-module $M_\l$ since $P_{n-k}$ is a projective
$\mS_{n-k}$-module. Then
$$\Ext_{\mS_n}^n(M_\l ,
\mS_n)\simeq P_{n-k}\t \mS_k / \sum_{i=1}^k(y_i-\l_i) \mS_k\simeq
P_{n-k}\t K\simeq P_{n-k}. $$ Then the projective dimension of the
$\mS_n$-module $M_\l$ is $k$. $\Box $

\begin{theorem}\label{28Dec8}
Let $M$ be a simple $\mS_n$-module, i.e. $M\simeq M_{\CN , \gm } =
P_{\CN } \t L_{C\CN } / \gm$ (Theorem \ref{17Dec8}). Then $\pd (M)
= |C\CN |$.
\end{theorem}

{\it Proof}. Without loss of generality we may assume that the
field $K=\bK$ is algebraically closed since $\pd_{\mS_n}(M)=
\pd_{\bK\t \mS_n}(\bK \t M)$. Then, up to order of the indices,
$M\simeq P_{n-k}\t (L_k/\sum_{i=1}^kL_k(y_i-\l_i))$ for some $k$
and a vector $\l = (\l_1, \ldots , \l_k)\in K^{*k}$. For each
$i=1,\ldots , k$, $K_{\l_i}:= L_1(i)/ L_1(i)(y_i-\l_i)$ is the
$\mS_1(i)$-module, and so $M\simeq P_{n-k}\t K_{(\l_1, \ldots ,
\l_k)}$ where $ K_{(\l_1, \ldots ,
\l_k)}:=\bigotimes_{i=1}^kK_{\l_i}$. We have to show that $\pd
(M)=k$. We deduce this fact from Corollary \ref{b10Jan9}. Let us
prove that, for all $s=0, 1,\ldots , k$, 
\begin{equation}\label{pdFKM}
\pd (P_{n-k}\t K_{(\l_1, \ldots , \l_s)}\t M_{(\l_{s+1}, \ldots ,
\l_k)})=k.
\end{equation}
To prove this fact we use induction on $s$. The case $s=0$ is
obvious, $\pd (P_{n-k}\t M_{(\l_1, \ldots , \l_k)})=k$ (Corollary
\ref{b10Jan9}). Suppose that $s>0$ and that (\ref{pdFKM}) holds
for all $s'<s$. Applying $-\t P_{n-k}\t K_{(\l_1, \ldots ,
\l_{s-1})}\t M_{(\l_{s+1}, \ldots , \l_k)}$ to the exact sequence
of $\mS_1(s)$-modules (see (\ref{MlKl}))
$$ 0\ra K[x_s]\ra M_{\l_s}\ra K_{\l_s}\ra 0$$
 we obtain the short exact sequence of $\mS_n$-modules:
\begin{eqnarray*}
  0& \ra & P_{n-k}\t K[x_s]\t K_{(\l_1, \ldots ,
\l_{s-1})}\t M_{(\l_{s+1}, \ldots , \l_k)}\ra
 P_{n-k}\t K_{(\l_1, \ldots ,
\l_{s-1})}\t M_{(\l_s, \ldots , \l_k)}\\
& \ra &  P_{n-k}\t K_{(\l_1, \ldots , \l_s)}\t M_{(\l_{s+1},
\ldots , \l_k)}\ra 0
\end{eqnarray*}
which can be written shortly as $0\ra X\ra Y\ra Z\ra 0$. Since
$\pd (X)\leq k-1$ and $\pd (Y)=k$ (by induction), we have $\pd
(Z)=k$. By induction, (\ref{pdFKM}) is true. In particular, $\pd
(M)=k$. This finishes the proof of the theorem. $\Box $

$\noindent $

\begin{corollary}\label{x28Dec8}
For all simple $\mS_n$-modules $M$, $\GK (M)+\pd (M) = \lgldim
(\mS_n)$.
\end{corollary}

{\it Proof}. Let $M$ be a simple $\mS_n$-module. Then $M\simeq
M_{\CN , \gm } $ (Theorem \ref{17Dec8}). By  Theorem
\ref{17Dec8}.(2,3), $\GK (M) = |\CN |$. By Theorem \ref{28Dec8},
$\pd (M) = |C\CN |$, and so $$\GK (M)+\pd (M) = |\CN |+|C\CN |=n=
\lgldim (\mS_n)\;\;  ({\rm Theorem}\;  \ref{A11Dec8}). \;\; \Box
$$

$\noindent $

\begin{lemma}\label{a26Dec8}
For each non-empty subset $\CN$ of the set $\{ 1, \ldots , n\}$
the intersection $\bigcap_{i\in \CN } \gp_i$ is a projective left
and right $\mS_n$-module. In particular, so are the ideals $\gp_1,
\ldots , \gp_n$.
\end{lemma}

{\it Proof}. We use induction on $n$ to prove the result. The case
$n=1$ is obvious since $F$ is a projective left and right
$\mS_1$-module (Proposition \ref{a11Dec8}.(2)). Suppose that $n>1$
and the result is true for all $n'<n$. The case $\CN = \{ 1,
\ldots ,n\}$ was considered already (Proposition
\ref{a11Dec8}.(2)). Without loss of generality we may assume that
$\CN = \{ 1, 2, \ldots , m\}$ and $m<n$. Then $\bigcap_{i=1 }^m
\gp_i =F_m\t \mS_{n-m}$ where $F_m$ is the projective left and
right $\mS_m$-module, by induction. Then it is obvious that $F_m\t
\mS_{n-m}$ is a  projective left and right $\mS_n$-module. $\Box $

$\noindent $

By Lemma \ref{a26Dec8}, for each number $s=1, \ldots , n$, we have
the projective $\mS_n$-module
$$I_s:= \bigoplus_{1\leq i_1<\cdots <i_s\leq n}\gp_{i_1}\cap \cdots
\cap \gp_{i_s}.$$ Consider the sequence of $\mS_n$-homomorphisms:
\begin{equation}\label{anres}
0\ra I_n\stackrel{d_n}{\ra} I_{n-1}\ra\cdots \ra
I_1\stackrel{d_1}{\ra} I_0:= \ga_n \ra 0
\end{equation}
where, for $s>1$,
$$d_s:\gp_{i_1}\cap \cdots \cap \gp_{i_s}\ra
\bigoplus_{t=1}^s \gp_{i_1}\cap \cdots \cap
\widehat{\gp_{i_t}}\cap \cdots \cap \gp_{i_s}, \;\; a\mapsto
((-1)^1a, (-1)^2a, \ldots, (-1)^sa), $$ and $d_1 (a_1, \ldots ,
a_n) = a_1+\cdots + a_n$ where $a_i\in \gp_i$ and the hat over a
symbol means that it is missed.

\begin{theorem}\label{29Dec8}
The sequence (\ref{anres}) is a projective resolution of the left
and right $\mS_n$-module $\ga_n$.
\end{theorem}

{\it Proof}. Since $\eta (I_s) = I_s$ and $\eta d_s= d_s \eta$ for
all $s$, it suffices to show that (\ref{anres}) is a projective
resolution of the left $\mS_n$-module $\ga_n$. By the very
definition of the maps $d_s$, $d_{s-1}d_s=0$ for all $s$. So, it
remains to prove the exactness of the complex (\ref{anres}). We
use induction on $n$. The case $n=2$ is obvious:
$$ 0\ra \gp_1\cap \gp_2\stackrel{d_2}{\ra} \gp_1\oplus
\gp_2\stackrel{d_1}{\ra}\gp_1+\gp_2\ra 0, \;\; d_2(a) = (-a,a),
\;\; d_1(u,v) = u+v.$$ So, let $n>2$, and we assume that the
result holds for all $n'<n$.

The idea of the proof is, first, to show  the exactness of the
complex at $s=1$, and then the exactness of the complex at $1<s<n$
 is deduced from the case $s=1$. Note that the complex is exact at
$s=n$ (since $d_n$ is an injection). For each $s$, $I_s=I_s'\oplus
I_s''$ where
$$I_s':= \bigoplus_{1\leq i_1<\cdots <i_s< n}\gp_{i_1}\cap \cdots
\cap \gp_{i_s}\;\; {\rm and}\;\; I_s'':= \bigoplus_{1\leq
i_1<\cdots <i_s= n}\gp_{i_1}\cap \cdots \cap \gp_{i_{s-1}}\cap
\gp_n.$$ Let us prove the exactness of the complex at $I_1$. By
induction, for $n-1$ the complex (\ref{anres}) is exact at
$I_1(n-1)$, i.e. the sequence of $\mS_{n-1}$-modules
\begin{equation}\label{aI2I1}
I_2(n-1)\stackrel{d_2(n-1)}{\ra}
I_1(n-1)\stackrel{d_1(n-1)}{\ra}\gp_1(n-1)+\cdots +\gp_{n-1}(n-1)
\end{equation}
is exact, i.e. $\ker ( d_1(n-1))=\im (d_2(n-1))$ where `$(n-1)$'
everywhere indicates that we consider the complex (\ref{anres})
for $n-1$. By applying $-\t\mS_1(n)$ to the sequence (\ref{aI2I1})
 we obtain the exact sequence 
\begin{equation}\label{I2I1}
I_2'\stackrel{d_2'}{\ra} I_1'\stackrel{d_1'}{\ra}\gp_1+\cdots
+\gp_{n-1}
\end{equation}
where $d_2'$ and $d_1'$ are the restrictions of the maps $d_2$ and
$d_1$ to $I_2'$ and $I_1'$ respectively. The sequence
$I_2\stackrel{d_2}{\ra} I_1\stackrel{d_1}{\ra}\ga_n$ can be
written as follows:
$$
I_2'\oplus I_2''\stackrel{\bigl(
\begin{smallmatrix}  d_2'& *\\ 0 & d_2''
\end{smallmatrix}\bigr)}{\ra} I_1'\oplus
\gp_n\stackrel{\bigl(
\begin{smallmatrix}  d_1'& {\rm id}\\
\end{smallmatrix}\bigr)}{\ra} \gp_1+\cdots +\gp_n
$$
where $I_2''=\bigoplus_{i=1}^{n-1}\gp_i\cap \gp_n$, $d_2''(a_1,
\ldots , a_{n-1})= -(a_1+\cdots + a_{n-1})$, $I_1'=
\gp_1\oplus\cdots \oplus \gp_{n-1}$. If an element $a= (a_1,
\ldots , a_n)\in I_1 = \gp_1\oplus \cdots \oplus \gp_n$ belongs to
$\ker (d_1)$, then
$$a_n = -a_1-\cdots - a_{n-1}\in (\gp_1+\cdots + \gp_{n-1})\cap
\gp_n= \gp_1\cap \gp_n+\cdots + \gp_{n-1}\cap \gp_n= \im
(d_2'').$$ and so $a_n= d_2''(b)$ for some element $b\in I_2''$.
Without loss of generality we may assume that $a_n=0$, i.e. $a\in
I_1'$ and $a\in \ker (d_1')$. Since the sequence (\ref{I2I1}) is
exact, we must have $a\in \im (d_2)$, and so the sequence
(\ref{anres}) is exact at $I_1$.

Now we use the (second) induction on $s$ to prove the exactness of
the sequence (\ref{anres}) at $I_s$. The case $s=1$ has been just
considered. So, let $s\geq 2$ and $s\neq n$, and we assume that
the complex (\ref{anres}) is exact at $I_{s'}$ for all
$s'=1,2,\ldots , s-1$. The sequence
$I_{s+1}\stackrel{d_{s+1}}{\ra} I_s\stackrel{d_s}{\ra}I_{s-1}$ can
be written as
$$
I_{s+1}'\oplus I_{s+1}''\stackrel{\bigl(
\begin{smallmatrix}  d_{s+1}'& *\\ 0 & d_{s+1}''
\end{smallmatrix}\bigr)}{\ra} I_s'\oplus
I_s''\stackrel{\bigl(
\begin{smallmatrix}  d_s'& *\\ 0 & d_s''
\end{smallmatrix}\bigr)}{\ra} I_{s-1}'\oplus I_{s-1}''
$$
where the maps $d_{s+1}'$ and $d_s'$ are the restrictions of the
maps $d_{s+1}$ and $d_s$ to $I_{s+1}'$ and $I_s'$ respectively. To
prove the exactness  at $I_s$ it suffices to prove the exactness
of two sequences
$I_{s+1}'\stackrel{d_{s+1}'}{\ra}I_s'\stackrel{d_s'}{\ra}I_{s-1}'$
and
 $I_{s+1}''\stackrel{d_{s+1}''}{\ra}I_s''\stackrel{d_s''}{\ra}I_{s-1}''$.
 By induction on $n$, the sequence
 $$I_{s+1}(n-1)\stackrel{d_{s+1}(n-1)}{\ra}I_s(n-1)\stackrel{d_s(n-1)}{\ra}I_{s-1}(n-1)$$
is exact. By applying $-\t\mS_1(n)$ to  this sequence   we obtain
the {\em exact} sequence
$I_{s+1}'\stackrel{d_{s+1}'}{\ra}I_s'\stackrel{d_s'}{\ra}I_{s-1}'$.

If $s=2$ then the sequence
$I_3''\stackrel{d_3''}{\ra}I_2''\stackrel{d_2''}{\ra}I_1''$ is
exact as it can be obtained by applying the exact functor
$\gp_n\t_{\mS_n}-$ to the exact sequence
$I_2'\stackrel{d_2'}{\ra}I_1'\stackrel{d_1'}{\ra}\mS_n$ (see
(\ref{I2I1})).

By induction on $s$, the sequence
$I_s\stackrel{d_s}{\ra}I_{s-1}\stackrel{d_{s-1}}{\ra}I_{s-2}$ is
exact. If $s>2$ then applying the exact functor
$-\t_{\mS_1(n)}F(n)$ to this sequence we obtain the {\em exact}
sequence
$I_{s+1}''\stackrel{d_{s+1}''}{\ra}I_s''\stackrel{d_s''}{\ra}I_{s-1}''$.
The functor $-\t_{\mS_1(n)}F(n)$ is exact since the
$\mS_1(n)$-module $F(n)$ is projective. The proof of the theorem
is complete. $\Box $


\section{Idempotent ideals of the algebra $\mS_n$}
\label{IDIDS}

In this section, all the idempotent ideals of the algebra $\mS_n$
are found (Theorem \ref{9Jan9}). It is proved that each idempotent
ideal distinct from $\mS_n$ is a unique product and a unique
intersection  of incomparable idempotent prime ideals (Theorem
\ref{27Ma7}). The intersection of idempotent ideals is always
equal to their product (Corollary \ref{a28Ma7}.(3)). The set
$\mI_n$ of all the idempotent ideals of the algebra $\mS_n$ is a
distributive lattice (Corollary \ref{28Ma7}).

Let $\CB_n$ be the set of all functions $f:\{ 1, 2, \ldots , n\}
\ra \F_2:= \{ 0,1\}$ where $\F_2:= \Z / 2\Z$ is a field. $\CB_n$
is a commutative ring with respect to addition and multiplication
of functions.  For $f,g\in \CB_n$, we write $f\geq g$ iff $f(i)
\geq g(i)$ for all $i=1, \ldots , n$ where $1>0$. Then $(\CB_n,
\geq )$ is a partially ordered set. For each function $f\in
\CB_n$, $I_f$ denotes the ideal $I_{f(1)}\t \cdots \t I_{f(n)}$ of
$\mS_n$ which is the tensor product of the ideals $I_{f(i)}$  of
the tensor components $\mS_1(i)$ in $\mS_n= \mS_1(1)\t \cdots \t
\mS_1(n)$ where $I_0:= F$ and $I_1:= \mS_1$. In particular,
$I_{(0,\ldots , 0)}=F_n$ and $I_{(1, \ldots , 1)}=\mS_n$.  $f\geq
g$ iff $I_f\supseteq I_g$. For $f,g\in \CB_n$, $I_fI_g= I_f\cap
I_g = I_{fg}$. Using induction on the number of functions we see
that, for $f_1, \ldots , f_s\in \CB_n$,
$$\prod_{i=1}^s I_{f_i}= \bigcap_{i=1}^s I_{f_i} = I_{f_1\cdots
f_s}.$$ Let $\CC_n$ be the set of all  subsets of $\CB_n$ all
distinct elements of which are incomparable (two distinct elements
$f$ and $g$ of $\CB_n$ are {\em incomparable} iff $f\not\leq g$
and $g\not\leq f$). For each $C\in \CC_n$, let $I_C:= \sum_{f\in
C}I_f$, the ideal of $\mS_n$.

Let  $\Sub_n$ be the set of all subsets of $\{ 1, \ldots, n \}$.
$\Sub_n$ is a partially ordered set with respect to `$\subseteq
$'. For each $f\in \CB_n$, the subset $\supp (f) := \{ i \, | \,
f(i) =1\}$ of $\{ 1, \ldots , n\}$ is called the {\em support} of
$f$. The map $ \CB_n\ra \Sub_n$, $ f\mapsto \supp (f)$,  is an
isomorphism of posets. Let $\SSub_n$ be the set of all subsets of
$\Sub_n$. An element $\{ X_1, \ldots , X_s\}$ of $\SSub_n$ is
called an {\em antichain} if for all $i\ne j$ such that $1\leq
i,j\leq s$ neither $X_i\subseteq X_j$ nor $X_i\supseteq X_j$. An
empty set and one element set are called antichains by definition.
Let $\Inc_n$ be the subset of $\SSub_n$ of all antichains  of
$\SSub_n$. Then the map
\begin{equation}\label{CnInc}
\CC_n\ra \Inc_n, \;\; \{ f_1, \ldots , f_s\} \mapsto \{ \supp
(f_1), \ldots , \supp (f_s)\},
\end{equation}
is a bijection.

$\noindent $

{\it Definition}. The number $\gd_n := |\Inc_n |$ is called the
{\em Dedekind} number.

$\noindent $

The Dedekind numbers  appeared in the paper of Dedekind
\cite{Dedekind-1871}. An asymptotic of the Dedekind numbers was
found by Korshunov \cite{Korshunov-1977}.

Let $\mI_n$ be the set of all the idempotent ideals of the algebra
$\mS_n$. The next theorem classifies all such ideals and gives a
canonical presentation for each of them.

\begin{theorem}\label{9Jan9}
\begin{enumerate}
\item The map $\CC_n \ra \mI_n$, $C\mapsto I_C:=\sum_{f\in C}I_f$,
is a bijection where $I_\emptyset :=0$. \item The set $\mI_n$ is
finite. Moreover, $|\mI_n |=\gd_n$ is the Dedekind number which
has the following asymptotic when $n\ra \infty$,
\cite{Korshunov-1977}:
$$
\gd_n\sim \begin{cases}
2^{n\choose \frac{n}{2}}e^{ {n\choose \frac{n}{2}-1} (2^{-\frac{n}{2}} +n^2 2^{-n-5}-n2^{-n-4})}&
 \text{if }n\;\; {\rm is \; even},\\
2\cdot 2^{n\choose \frac{n-1}{2}}e^{\frac{n(n-3)}{2}
(2^{-\frac{n+3}{2}} +n^2 2^{-n-6})    + \frac{n(n-1)}{2} (
2^{-\frac{n+1}{2}} +n^22^{-n-4})}& \text{if }n\;\; {\rm is \;
odd}.
\end{cases}
$$
\item $\eta (I) = I$ for all idempotent  ideals $I$
 of the algebra $\mS_n$.
\end{enumerate}
\end{theorem}

{\it Proof}. 1. The map $C\mapsto I_C$ is well defined since
$$I^2= (\sum_{f\in C}I_f)^2= \sum_{f\in C}I^2_f+\sum_{f\neq g} I_fI_g=\sum_{f\in
C}I_f+\sum_{f\neq g} I_fI_g=\sum_{f\in C}I_f=I.$$ Choose a basis
for the algebra $\mS_n=\t_{i=1}^n\mS_1(i)$ which is the tensor
product of bases of the tensor components $\mS_1(i)$, and each
basis for $\mS_1(i)$ is an extension of a basis for its  subspace
$F(i) \subseteq \mS_1(i)$. Then it is obvious that the map
$C\mapsto I_C$ is injective.

Clearly, the algebra $\mS_n$ is an idempotent ideal. So, it
remains to show that if $I$ is an idempotent ideal of the algebra
$\mS_n$ such that $I\neq \mS_n$ then $I=I_C$ for some $C$.

First, let us show that $I\subseteq \ga_n$. The image $\bI$ of the
ideal $I$ under the epimorphisms $\mS_n\ra \mS_n/\ga_n\simeq L_n$
is an idempotent ideal, hence either $\bI =0$ or $\bI = L_n$. The
second case is impossible, otherwise we would have
$I+\ga_n=\mS_n$. Let $\gm$ be a maximal ideal of $\mS_n$ that
contains $I$ (it exists since $I\neq \mS_n$). Then
$\mS_n=I+\ga_n\subseteq \gm$ (Corollary \ref{a9Jan9}), a
contradiction. Therefore, $I\subseteq \ga_n$.

Now, the result is obvious for $n=1$ (Proposition
\ref{10Dec8}.(3)). So, let $n>1$ and we assume that the result is
true for all $n'<n$. Multiplying the chain of  inclusions
$I\subseteq \ga_n\subseteq \mS_n$ on the left by the ideal $I$ we
have $I=I^2\subseteq I\ga_n\subseteq I\mS_n = I$, and so
$$ I=I\ga_n = I(\gp_1+\cdots +\gp_n ) = I\gp_1+\cdots + I\gp_n.$$
By Lemma \ref{c23Dec8}.(1), $I\cap \gp_1 = I\cap (F\t \mS_{n-1}) =
F\t I_{n-1}$ for an ideal $I_{n-1}$ of the algebra $\mS_{n-1}$.
Clearly,
$$ \gp_1I\subseteq I\cap \gp_1\subseteq F\t I_{n-1}\subseteq \gp_1
I,$$ and so $I\cap \gp_1= \gp_1I$.  By symmetry, $I\cap \gp_i =
\gp_i I$ for all $i$. Now, $I= \sum_{i=1}^n I\gp_i= \sum_{i=1}^n
I\cap \gp_i$. To finish the proof it suffices to show that each
ideal $I\cap \gp_i$ has the form $I_{C_i}$ for some $C_i$. By
symmetry, it suffices to prove this for $i=1$. Recall that the
ideals of the algebra $\mS_n$ commute (Theorem \ref{b23Dec8}). On
the one hand, $(I\cap \gp_1)^2= (I\gp_1)^2= I^2\gp_1^2= I\gp_1=
I\cap \gp_1= F\t I_{n-1}$. On the other hand, $(I\cap
\gp_1)^2=(F\t I_{n-1})^2= F^2\t I_{n-1}^2= F\t I^2_{n-1}$.
Therefore, $I_{n-1}$ is an idempotent ideal of the algebra
$\mS_{n-1}$. By induction, the ideal $I_{n-1}$ has the  required
form, and so the ideal $I\cap \gp_1$ has the form $I_{C_1}$. This
finishes the proof of statement 1.

2. Statement 2 follows from statement 1 and (\ref{CnInc}).

3. By statement 1, $I=I_C$ for some $S$. Then $\eta (I_C)= \eta
(\sum_{f\in C}I_f ) =\sum_{f\in C}\eta (I_f) =\sum_{f\in C} I_f=
I_C$ since $\eta (I_f) = I_f$ for all $f$. $\Box $

$\noindent $

$(\Spec (\mS_n) , \subseteq )$ is a poset. Two primes $\gp $ and
$\gq$ are called {\em incomparable} if neither $\gp \subseteq \gq
$ nor $\gp \supseteq \gq$.

For each idempotent ideal $\ga$ of $\mS_n$ such that $\ga\neq
\mS_n$, let $\Min (\ga)$ be the set of all minimal primes over
$\ga$. The set $\Min (\ga ) $ is a non-empty finite set (Theorem
\ref{25Feb9}) and each element of $\Min (\ga )$ is an {\em
idempotent}, prime ideal (Theorem \ref{27Ma7}). The proof of
Theorem \ref{27Ma7} provides a direct, short proof of the fact
that  the set $\Min (\ga )$ is finite and non-empty for an each
{\em idempotent} ideal $\ga$ of $\mS_n$.

For each $f\in \CB_n$, the set $\csupp (f):= \{ i\, | \, f(i)=0\}$
is called the {\em co-support} of $f$. Clearly, $\csupp (f) = \{
1, \ldots , n\} \backslash \supp (f)$.

\begin{theorem}\label{27Ma7}
\begin{enumerate}
\item Each idempotent ideal $\ga$ of $\mS_n$ such that $\ga \neq
\mS_n$ is a unique product of incomparable idempotent primes, i.e.
if $\ga = \gq_1\cdots \gq_s= \gr_1\cdots \gr_t$ are two such
products then $s=t$ and $\gq_1= \gr_{\s (1)}, \ldots , \gq_s=
\gr_{\s (s)}$ for a permutation $\s$ of $\{ 1, \ldots, n\}$. \item
Each idempotent ideal $\ga$ of $\mS_n$ such that $\ga \neq \mS_n$
is a unique intersection of incomparable idempotent primes, i.e.
if $\ga = \gq_1\cap \cdots\cap  \gq_s= \gr_1\cap \cdots \cap
\gr_t$ are two such intersections then $s=t$ and $\gq_1= \gr_{\s
(1)}, \ldots , \gq_s= \gr_{\s (s)}$ for a permutation $\s$ of $\{
1, \ldots, n\}$. \item For each idempotent ideal $\ga$ of $\mS_n$
such that $\ga \neq \mS_n$, the sets of incomparable idempotent
primes in statements 1 and 2 are the same, and so $\ga=\gq_1\cdots
\gq_s=\gq_1\cap \cdots\cap \gq_s$. \item The ideals $\gq_1,\ldots
, \gq_s$ in statement 3 are the minimal  primes of $\ga$, and so
$\ga = \prod_{\gp \in \Min (\ga )}\gp =\cap_{\gp \in \Min (\ga
)}\gp$. In particular, each element of $\Min (\ga )$ is an
idempotent prime ideal of $\mS_n$.
\end{enumerate}
\end{theorem}

{\it Proof}. 1. For each idempotent ideal $\ga$ of $\mS_n$ such
that $\ga \neq \mS_n$, we have to prove that $\ga$ is a product of
incomparable idempotent primes and that this product is unique.
Since the ring $\mS_n$ is prime these two statements are obvious
when $\ga =0$. So, let $\ga \neq 0$.

{\em Existence}: Let $f\in \CB_n$; then $I_f=\prod_{i\in \csupp
(f)}\gp_i$. Let $\gb$ be any idempotent ideal of $\mS_n$. Since
$\gb^2= \gb$,
it follows at once that 
\begin{equation}\label{Ifb}
I_f+\gb = \prod_{i\in \csupp (f)} (\gp_i+\gb ).
\end{equation}
By Theorem \ref{9Jan9}.(1),  $\ga = I_{f_1} +\cdots + I_{f_s}$ for
some $f_i\in \CB_n$. Repeating $s$ times (\ref{Ifb}), we see that
\begin{equation}\label{aiSup}
\ga = \prod_{i_1\in \csupp (f_1),\ldots , i_s\in \csupp
(f_s)}(\gp_{i_1} +\cdots + \gp_{i_s})
\end{equation}
is the product of idempotent primes, by Corollary
\ref{g21Dec8}.(1,2). It follows that $\Min (\ga ) $ is a non-empty
finite set each element of which is an idempotent prime ideal.
Note that the ideals of $\mS_n$ commute, and if $\gp \subseteq
\gq$ is an inclusion of idempotent primes then $\gp \gq = \gp$
(Corollary \ref{g21Dec8}.(4)). Using these facts and
(\ref{aiSup}), we see that $\ga$ is a product of {\em
incomparable} idempotent primes.

Uniqueness follows from the next lemma which will be used several
times in the proof of this theorem.

\begin{lemma}\label{u28Ma7}
Let $\{ \gq_1, \ldots , \gq_s\}$ and $\{ \gr_1, \ldots , \gr_t\}$
be two sets of incomparable ideals of a ring such that each ideal
from the first set contains an ideal from the second and each
ideal from the second set contains an ideal from the first. Then
$s=t$ and $\gq_1= \gr_{\s (1)}, \ldots , \gq_s= \gr_{\s (s)}$ for
a permutation $\s$ of $\{ 1, \ldots, n\}$.
\end{lemma}

{\it Proof of Lemma \ref{u28Ma7}}. For each $\gq_i$, there are
ideals $\gr_j$ and $\gr_k$ such that $\gr_j\subseteq
\gq_i\subseteq \gr_k$, hence $\gq_i= \gr_j = \gr_k$ since the
ideals $\gr_j$ and $\gr_k$ are incomparable if distinct. This
proves that for each ideal $\gq_i$ there exists a unique ideal,
say $\gr_{\s (i)}$, such that $\gq_i= \gr_{\s (i)}$. By symmetry,
for each ideal $\gr_j$ there exists a unique  ideal, say
$\gq_{\tau (j)}$, such that $\gr_j = \gq_{\tau (j)}$. Then, $s=t$
and $\gq_1= \gr_{\s (1)}, \ldots , \gq_s= \gr_{\s (s)}$ for the
permutation $\s$ of $\{ 1, \ldots, n\}$. $\Box $

{\em Uniqueness}: Let $\ga = \gq_1\cdots \gq_s= \gr_1\cdots \gr_t$
 be  two products of  incomparable idempotent primes. Each ideal $\gq_i$ contains
an   ideal $\gr_j$, and each ideal $\gr_k$ contains an ideal
$\gq_l$. By Lemma \ref{u28Ma7}, $s=t$ and $\gq_1= \gr_{\s (1)},
\ldots , \gq_s= \gr_{\s (s)}$ for a permutation $\s$ of $\{ 1,
\ldots, s\}$.

2. {\em Uniqueness}: Suppose that an ideal $\ga$ has two
presentations $\ga = \gq_1\cap \cdots \cap\gq_s= \gr_1\cap
\cdots\cap  \gr_t$ of incomparable idempotent primes. The sets $\{
\gq_1, \ldots , \gq_s\}$ and $\{ \gr_1, \ldots , \gr_t\}$ of
incomparable  idempotent primes satisfy the conditions of Lemma
\ref{u28Ma7}, and so uniqueness follows.

{\em Existence}: Let $\mI'$ be the set of idempotent ideals of
$\mS_n$ that are intersection of incomparable idempotent primes.
Then $\mI'\subseteq \mI_n$. The map
$$\mI_n \ra \mI',\;\; \gq_1\cdots \gq_s\mapsto \gq_1\cap \cdots \cap
\gq_s,$$ is a  bijection since $|\mI_n |<\infty$ and by uniqueness
of presentations $\gq_1\cdots \gq_s$ (statement 1) and $\gq_1\cap
\cdots \cap \gq_s$ (see above) where $\gq_1, \ldots , \gq_s$ are
incomparable idempotent primes. Then $\mI_n = \mI'$. This proves
that each idempotent ideal $\ga$ of $\mS_n$ is  an intersection of
incomparable idempotent primes.

3. Let $\ga$ be an idempotent ideal of $\mS_n$ and $\ga =
\gq_1\cdots \gq_s= \gr_1\cap \cdots \cap \gr_t$ where $S:=\{
\gq_1, \ldots , \gq_s\}$ and $T:=\{ \gr_1, \ldots , \gr_t\}$ are
sets of incomparable idempotent primes. The sets $S$ and $T$
satisfy the conditions of Lemma \ref{u28Ma7}, and so $s=t$ and
$\gq_1= \gr_{\s (1)}, \ldots , \gq_s= \gr_{\s (s)}$ for a
permutation $\s$ of $\{ 1, \ldots, n\}$. This means that $\ga =
\gq_1\cdots \gq_s= \gq_1\cap \cdots \cap \gq_s$.

4. Let $\ga = \gq_1\cdots \gq_s = \gq_1\cap \cdots \cap \gq_s$ be
as in statement 3 and let $\Min (\ga )=\{ \gr_1, \ldots , \gr_t\}$
be the set of minimal primes over $\ga$. Then $\Min (\ga )
\subseteq S:=\{ \gq_1, \ldots \gq_s\}$ ($\ga = \gq_1\cdots
\gq_s\subseteq \gr_i$ implies $\gq_j\subseteq \gr_i$ for some $j$,
and so $\gq_j = \gr_i$ by the minimality of $\gr_i$). Up to order,
let $\gr_1= \gq_1, \ldots , \gr_t = \gq_t$. It remains to show
that $t=s$. Suppose that $t<s$, we seek a contradiction. This
means that each idempotent prime $\gq_i$, $i=t+1, \ldots , s$,
contains $\ga$ and is {\em not} a minimal prime  over $\ga$.
Hence, $\gq_i$ contains a minimal idempotent prime, say $\gq_{\tau
(i)}$, a contradiction (the ideals $\gq_i$ and $\gq_{\tau (i)}$
are incomparable).  $\Box $

\begin{corollary}\label{a28Ma7}
Let  $\ga$ and $\gb$ be idempotent ideals of $\mS_n$ distinct from
$\mS_n$ in statement 1, 2 and 5. Then
\begin{enumerate}
\item $\ga = \gb$ iff $\Min (\ga ) = \Min (\gb )$.\item $\Min (\ga
\cap \gb ) = \Min (\ga \gb ) = $ the set of minimal elements (with
respect to inclusion) of the set $\Min (\ga ) \cup \Min (\gb )$.
\item $ \ga\cap \gb = \ga \gb$. \item If $\ga \subseteq \gb$ then
$\ga \gb = \ga$. \item $\ga \subseteq \gb$ iff $\Min (\ga )
 \eqslantless \Min (\gb )$ (the $\eqslantless$ means that
and each $\gq\in \Min (\gb )$ contains some $\gp\in \Min (\ga )$).
\end{enumerate}
\end{corollary}

{\it Proof}. 1. Statement 1 is obvious due to Theorem
\ref{27Ma7}.(4).

2. Let $\CM$ be the set of minimal elements of the union $\Min
(\ga ) \cup \Min (\gb )$. The elements of $\CM$ are incomparable,
and (by Theorem \ref{27Ma7}.(4))
$$ \ga \cap \gb = \cap_{\gp\in \Min (\ga )}\cap \cap_{\gq \in \Min
(\gb )}\gq = \cap_{\gr \in \CM } \gr.$$ By Theorem
\ref{27Ma7}.(2), $\Min (\ga \cap \gb ) = \CM$. By Corollary
\ref{g21Dec8}.(4), $$\ga \gb = \prod_{\gp \in \Min (\ga ) }\gp
\cdot \prod_{\gq \in \Min (\gb )}\gq = \prod_{\gr \in \CM}\gr =
\ga \cap \gb.$$

3. The result is obvious if one of  the ideals is equal to
$\mS_n$. So, let the ideals are distinct from $\mS_n$.  By
statement 2, $\Min (\ga \cap \gb ) = \Min (\ga \gb )$, then, by
statement 1, $\ga\cap \gb = \ga \gb$.

4. If $\ga \subseteq \gb$ then, by statement 3, $\ga \gb = \ga\cap
\gb = \ga$.

5. $(\Rightarrow )$ If $\ga \subseteq \gb$ then $\Min (\ga )
\eqslantless\Min (\gb )$ since $\ga = \prod_{\gp\in \Min (\ga
)}\gp \subseteq \prod_{\gq \in \Min (\gb )}\gq = \gb$.

$(\Leftarrow )$ Suppose that $\Min (\ga ) \eqslantless\Min (\gb
)$. For each $\gq\in \Min (\gb )$, let $S(\gq )$ be the set
(necessarily non-empty) of $\gp\in \Min (\ga )$ such that $\gp
\subseteq \gq$. Then $\Min (\ga ) \supseteq S:= \cup_{\gq\in \Min
(\gb )} S(\gq )$ and
$$\ga =\cap_{\gp \in \Min (\ga )}\gp \subseteq \cap_{\gp \in S}\gp \subseteq
\cap_{\gq\in \Min (\gb )}\gq = \gb.\;\;\; \Box $$

\begin{corollary}\label{28Ma7}
The lattice $\mI_n$  of idempotent ideals of the algebra $\mS_n$
is distributive, i.e. $(\ga \cap \gb ) \gc = \ga \gc\cap \gb \gc$
for all ideals $\ga$, $\gb$, and $\gc$.
\end{corollary}

{\it Proof}.  By Corollary \ref{a28Ma7}.(3), $(\ga \cap \gb ) \gc
= \ga \cap \gb \cap \gc = (\ga \cap \gc ) \cap ( \gb \cap \gc ) =
\ga \gc \cap \gb \gc$.  $\Box $

\begin{theorem}\label{1Jun7}
Let  $\ga$ be an  idempotent ideal of $\mS_n$, and $\CM$ be the
set of  minimal elements with respect to inclusion of the  set of
minimal primes of idempotent ideals $\ga_1, \ldots , \ga_k$ of
$\mS_n$. Then
\begin{enumerate}
\item  $\ga = \ga_1\cdots \ga_k$ iff $\Min (\ga ) = \CM$.\item
$\ga = \ga_1\cap \cdots\cap \ga_k$ iff $\Min (\ga ) = \CM$.
\end{enumerate}
\end{theorem}

{\it Proof}. By Corollary \ref{a28Ma7}.(3), it suffices to prove,
say, the first statement.

$(\Rightarrow )$ Suppose that $\ga = \ga_1\cdots \ga_k$ then, by
Theorem \ref{27Ma7}.(4) and Corollary \ref{a28Ma7}.(4), $$\ga =
\prod_{i=1}^k\prod_{\gq_{ij}\in \Min (\ga_i)} \gq_{ij} =
\prod_{\gq \in \CM} \gq,$$
 and so $\Min (\ga ) = \CM$, by Theorem
\ref{27Ma7}.(4).

$(\Leftarrow )$ If $\Min (\ga ) = \CM$ then, by Corollary
\ref{a28Ma7}.(4), $\ga = \ga_1\cdots \ga_k$.  $\Box $

$${\bf Acknowledgements}$$

The author would like to thank the referee for comments.

Department of Pure Mathematics

University of Sheffield

Hicks Building

Sheffield S3 7RH

UK

email: v.bavula@sheffield.ac.uk

\end{document}